%This is a plain tex file with the paper
% {On the Very Weak 0-1 Law for Random Graphs with Orders}
% by Saharon SHELAH

\ifx\shlhetal\undefinedcontrolsequence\let\shlhetal\relax\fi

%=============================================================================
%
%    LNM.STY  TEX-Makros fr "Lecture Notes in Mathematics" (22.Dez.1989)
%
%=============================================================================
%
\hsize=144mm
\newdimen\seitenbreite
\seitenbreite=\hsize
\topskip=20mm
\vsize=230mm
\frenchspacing
\raggedbottom
\parskip=0pt
\baselineskip=6mm
\lineskip=1pt %0pt
\lineskiplimit=0pt %-100pt
\overfullrule=3pt
%\nopagenumbers
\pretolerance=200
\tolerance=500
\widowpenalty=10000
\clubpenalty=10000
%
%---------- Fonts ------------------------------------------------------------
%
\font\small=cmr10 
\font\tanorm=cmr12 scaled \magstep3
\font\taeins=cmmi10 scaled \magstep4
\font\tazwei=cmsy10 scaled \magstep4
%-------------------------------------------
\font\tbnorm=cmbx12 scaled \magstep1
\font\tbeins=cmmib10 scaled \magstep2
\font\tbzwei=cmbsy10 scaled \magstep2
%-------------------------------------------
\font\norm=cmr12
\font\eins=cmmi12
\font\zwei=cmsy10 scaled 1200
\font\drei=cmex10 scaled 1200
\font\msxm=msxm10 scaled 1200
\font\msym=msym10 scaled 1200
%-------------------------------------------
\font\klnull=cmr10
\font\kleins=cmmi10
\font\klzwei=cmsy10
\font\kldrei=cmex10
%-------------------------------------------
\font\snull=cmr9
\font\seins=cmmi9
\font\szwei=cmsy9
\font\sdrei=cmex10 scaled 900
%-------------------------------------------
\font\ssnull=cmr7
\font\sseins=cmmi7
\font\sszwei=cmsy7
\font\ssdrei=cmex10 scaled 700
%-------------------------------------------
\font\klbf=cmb10
\font\klku=cmti10
\font\ku=cmti12
\let\it=\ku

\let\bf=\schmalbf

\norm
\textfont0=\norm
\textfont1=\eins
\textfont2=\zwei
\textfont3=\drei
\scriptfont0=\snull
\scriptfont1=\seins
\scriptfont2=\szwei
\scriptfont3=\sdrei
\scriptscriptfont0=\ssnull
\scriptscriptfont1=\sseins
\scriptscriptfont2=\sszwei
\scriptscriptfont3=\ssdrei
%
%---------- Kolumnentitel ----------------------------------------------------
%
\newcount\anfseite
\def\dummyheadline{\hss}

\def\leftheadline{\klnull\textfont1=\kleins\textfont2=\klzwei
 \ifnum\pageno<0 \uppercase\romannumeral-\pageno \else\number\pageno \fi
 \hskip9mm\lkol\hfill}
\def\rightheadline{\klnull\textfont1=\kleins\textfont2=\klzwei\hfill\rkol\hskip9mm
 \ifnum\pageno<0 \uppercase\romannumeral-\pageno \else\number\pageno \fi}
\headline={\vbox{\vglue-1mm\noindent\ifodd\pageno
 \ifnum\pageno=\anfseite\dummyheadline\else\rightheadline\fi
 \else
 \ifnum\pageno=\anfseite\dummyheadline\else\leftheadline\fi
 \fi\vskip5mm~}}
\def\lkol{}
\def\rkol{}
%
%---------- Kapitel 1. Ordnung -----------------------------------------------
%
\def\ta#1#2{\par\noindent\vbox to 50mm{\bgroup\par\noindent
 \baselineskip=24pt
 \tanorm
 \textfont0=\tanorm
 \textfont1=\taeins
 \textfont2=\tazwei
 \raggedright #1~~#2 \vfill\egroup}}
%
%---------- Kapitel 2. Ordnung -----------------------------------------------
%
\def\tb#1#2{\vskip28.8pt\par\noindent\bgroup\baselineskip=17.3pt
 \tbnorm
 \textfont0=\tbnorm
 \textfont1=\tbeins
 \textfont2=\tbzwei
 \raggedright
 #1~~#2 \vskip14.4pt
 \par\noindent\egroup\ignorespaces}
%
%---------- Fuanote ----------------------------------------------------------
%
\def\fn#1#2{\bgroup
	  \baselineskip=12pt
	  \lineskip=0pt
	  \lineskiplimit=0pt
	  \klnull
	  \textfont0=\klnull
	  \textfont1=\kleins
	  \textfont2=\klzwei
	  \textfont3=\kldrei
	  \scriptfont0=\ssnull
	  \let\ku=\klku
	  \let\bf=\klbf
	  \footnote{$^{#1}$}{#2}\egroup}
\def\footnoterule{\kern-3pt\hrule width 20mm \kern 2.6pt}
%
%---------- Kleindruck -------------------------------------------------------
%
\def\akl{\bgroup
	    \par
	    \noindent
	    \baselineskip=12pt
	    \klnull
	    \textfont0=\klnull
	    \textfont1=\kleins
	    \textfont2=\klzwei
	    \textfont3=\kldrei
	    \scriptfont0=\ssnull
	    \let\bf=\klbf
	    \let\ku=\klku}
\def\ekl{\par
	    \noindent
	    \egroup}
%
%---------- Literatur --------------------------------------------------------
%
\def\ref{\par
	 \noindent
	 \hangindent=4mm
	 \hangafter=1
	 \baselineskip=12pt
	  \klnull
	    \textfont0=\klnull
	    \textfont1=\kleins
	    \textfont2=\klzwei
	    \textfont3=\kldrei
	    \scriptfont0=\ssnull
	 \let\bf=\klbf
	 \let\ku=\klku}
%
%-----------------------------------------------------------------------------
%

			    % Absatz
	    % neue Seite beginnen
\def\lz{\vskip6mm\noindent}	    % Leerzeile
\def\hz{\vskip7.2pt\noindent}	    % halbe Leerzeile
     % dreiviertel Leerzeile
	    % viertel Leerzeile
\def\nz{\hfil\break\noindent}	    % neue Zeile beginnen
\def\ts{\thinspace}		    % kleiner horizontaler Zwischenraum
		    % Leerstelle fr eine Ziffer
		    % Leerstelle fr zwei Ziffern
%
%-----------------------------------------------------------------------------
%

%
\def\lq{\hbox{\norm``}}
\def\rq{\hbox{\norm''}}

\def\mod{\hbox{\norm mod}}

\def\Rang{\hbox{\norm Rang}}

\def\lk{\langle}
\def\rk{\rangle}

\def\conc{\widehat{~~}}
\def\uhr{\hbox{\msxm\char'026}}

\def\qed{\hfill\hbox{\msxm\char'003}}
\def\Vdash{\hbox{\ts\msxm\char'015}}

\def\paragraph{\S}
\def\atil{a\llap{\lower4pt\hbox{$\textfont2=\sszwei\sim$}}{}}
\def\ftil{f\llap{\lower6pt\hbox{$\textfont2=\sszwei\sim$\kern1pt}}{}}
\def\gtil{g\llap{\lower6pt\hbox{$\textfont2=\sszwei\sim$}}{}}
\def\htil{h\llap{\lower6pt\hbox{$\textfont2=\sszwei\sim$}}{}}
\def\ktil{k\llap{\lower6pt\hbox{$\textfont2=\sszwei\sim$}}{}}
\def\ntil{n\llap{\lower4pt\hbox{$\textfont2=\sszwei\sim$}}{}}
\def\qtil{q\llap{\lower6pt\hbox{$\textfont2=\sszwei\sim$}}{}}
\def\rtil{r\llap{\lower4pt\hbox{$\textfont2=\sszwei\sim$}}{}}
\def\xtil{x\llap{\lower4pt\hbox{$\textfont2=\sszwei\sim$}}{}}
\def\ytil{y\llap{\lower6pt\hbox{$\textfont2=\sszwei\sim$}}{}}
\def\Atil{A\llap{\lower4pt\hbox{$\textfont2=\sszwei\sim$\kern2pt}}{}}
\def\Btil{B\llap{\lower4pt\hbox{$\textfont2=\sszwei\sim$\kern2pt}}{}}
\def\Ctil{C\llap{\lower4pt\hbox{$\textfont2=\sszwei\sim$\kern2pt}}{}}
\def\Dtil{D\llap{\lower4pt\hbox{$\textfont2=\sszwei\sim$\kern2pt}}{}}
\def\Etil{E\llap{\lower4pt\hbox{$\textfont2=\sszwei\sim$\kern2pt}}{}}
\def\Ftil{F\llap{\lower4pt\hbox{$\textfont2=\sszwei\sim$\kern2pt}}{}}
\def\Gtil{G\llap{\lower4pt\hbox{$\textfont2=\sszwei\sim$\kern2pt}}{}}
\def\Mtil{M\llap{\lower4pt\hbox{$\textfont2=\sszwei\sim$\kern2pt}}{}}
\def\Ptil{P\llap{\lower4pt\hbox{$\textfont2=\sszwei\sim$\kern2pt}}{}}
\def\Qtil{Q\llap{\lower5pt\hbox{$\textfont2=\sszwei\sim$\kern2pt}}{}}
\def\Rtil{R\llap{\lower4pt\hbox{$\textfont2=\sszwei\sim$\kern2pt}}{}}
\def\Stil{S\llap{\lower4pt\hbox{$\textfont2=\sszwei\sim$\kern2pt}}{}}
\def\Ttil{T\llap{\lower4pt\hbox{$\textfont2=\sszwei\sim$\kern2pt}}{}}
\def\Xtil{X\llap{\lower4pt\hbox{$\textfont2=\sszwei\sim$\kern2pt}}{}}
\def\Ytil{Y\llap{\lower4pt\hbox{$\textfont2=\sszwei\sim$\kern2pt}}{}}
\def\alphatil{\alpha\llap{\lower4pt\hbox{$\textfont2=\sszwei\sim$}}{}}
\def\betatil{\beta\llap{\lower6pt\hbox{$\textfont2=\sszwei\sim$}}{}}
\def\kappatil{\kappa\llap{\lower4pt\hbox{$\textfont2=\sszwei\sim$}}{}}
\def\lambdatil{\lambda\llap{\lower4pt\hbox{$\textfont2=\sszwei\sim$}}{}}
\def\mutil{\mu\llap{\lower6pt\hbox{$\textfont2=\sszwei\sim$}}{}}
\def\phitil{\phi\llap{\lower6pt\hbox{$\textfont2=\sszwei\sim$}}{}}
\def\psitil{\psi\llap{\lower6pt\hbox{$\textfont2=\sszwei\sim$}}{}}
\def\tautil{\tau\llap{\lower4pt\hbox{$\textfont2=\sszwei\sim$}}{}}
\def\greatil{>\llap{\lower4pt\hbox{$\textfont2=\sszwei\sim$}}{}}
\def\membertil{\in\llap{\lower4pt\hbox{$\textfont2=\sszwei\sim$}}{}}
\def\SDtil{{\cal D}\llap{\lower4pt\hbox{$\textfont2=\sszwei\sim$}}{}}
\def\SA{{\cal A}}

\def\SI{{\cal I}}

\def\SL{{\cal L}}

%
%-----------------------------------------------------------------------------
%
\abovedisplayskip=15pt
\belowdisplayskip=15pt
%
%-----------------------------------------------------------------------------
%
%\catcode`\@=11
% nur wichtig, wenn Kommandos @ enthalten
%\catcode`\@=12
\noindent
%=============================================================================

% Newest version!!
\def\today{\ifcase\month\or January\or February\or March\or April\or
    May\or June\or July\or August\or September\or October\or November
    \or December\fi\space\number\day, \number\year}

%capital letters with a tilde under them, and a subset.

\def\Atill_#1{\setbox5=\hbox{$A_#1$}\copy5\kern-\wd5%
\lower4pt\hbox to \wd5 {$\textfont2=\sszwei\sim$\hfil}}
\def\Btill_#1{\setbox5=\hbox{$B_#1$}\copy5\kern-\wd5%
\lower4pt\hbox to \wd5 {$\textfont2=\sszwei\sim$\hfil}}
\def\Ctill_#1{\setbox5=\hbox{$C_#1$}\copy5\kern-\wd5%
\lower4pt\hbox to \wd5 {$\textfont2=\sszwei\sim$\hfil}}
\def\Dtill_#1{\setbox5=\hbox{$D_#1$}\copy5\kern-\wd5%
\lower4pt\hbox to \wd5 {$\textfont2=\sszwei\sim$\hfil}}
\def\Etill_#1{\setbox5=\hbox{$E_#1$}\copy5\kern-\wd5%
\lower4pt\hbox to \wd5 {$\textfont2=\sszwei\sim$\hfil}}
\def\Ftill_#1{\setbox5=\hbox{$F_#1$}\copy5\kern-\wd5%
\lower4pt\hbox to \wd5 {$\textfont2=\sszwei\sim$\hfil}}
\def\Gtill_#1{\setbox5=\hbox{$G_#1$}\copy5\kern-\wd5%
\lower4pt\hbox to \wd5 {$\textfont2=\sszwei\sim$\hfil}}
\def\Htill_#1{\setbox5=\hbox{$H_#1$}\copy5\kern-\wd5%
\lower4pt\hbox to \wd5 {$\textfont2=\sszwei\sim$\hfil}}
\def\Itill_#1{\setbox5=\hbox{$I_#1$}\copy5\kern-\wd5%
\lower4pt\hbox to \wd5 {$\textfont2=\sszwei\sim$\hfil}}
\def\Jtill_#1{\setbox5=\hbox{$J_#1$}\copy5\kern-\wd5%
\lower4pt\hbox to \wd5 {$\textfont2=\sszwei\sim$\hfil}}
\def\Ktill_#1{\setbox5=\hbox{$K_#1$}\copy5\kern-\wd5%
\lower4pt\hbox to \wd5 {$\textfont2=\sszwei\sim$\hfil}}
\def\Ltill_#1{\setbox5=\hbox{$L_#1$}\copy5\kern-\wd5%
\lower4pt\hbox to \wd5 {$\textfont2=\sszwei\sim$\hfil}}
\def\Mtill_#1{\setbox5=\hbox{$M_#1$}\copy5\kern-\wd5%
\lower4pt\hbox to \wd5 {$\textfont2=\sszwei\sim$\hfil}}
\def\Ntill_#1{\setbox5=\hbox{$N_#1$}\copy5\kern-\wd5%
\lower4pt\hbox to \wd5 {$\textfont2=\sszwei\sim$\hfil}}
\def\Otill_#1{\setbox5=\hbox{$O_#1$}\copy5\kern-\wd5%
\lower4pt\hbox to \wd5 {$\textfont2=\sszwei\sim$\hfil}}
\def\Ptill_#1{\setbox5=\hbox{$P_#1$}\copy5\kern-\wd5%
\lower4pt\hbox to \wd5 {$\textfont2=\sszwei\sim$\hfil}}
\def\Qtill_#1{\setbox5=\hbox{$Q_#1$}\copy5\kern-\wd5%
\lower4pt\hbox to \wd5 {$\textfont2=\sszwei\sim$\hfil}}
\def\Rtill_#1{\setbox5=\hbox{$R_#1$}\copy5\kern-\wd5%
\lower4pt\hbox to \wd5 {$\textfont2=\sszwei\sim$\hfil}}
\def\Still_#1{\setbox5=\hbox{$S_#1$}\copy5\kern-\wd5%
\lower4pt\hbox to \wd5 {$\textfont2=\sszwei\sim$\hfil}}
\def\Ttill_#1{\setbox5=\hbox{$T_#1$}\copy5\kern-\wd5%
\lower4pt\hbox to \wd5 {$\textfont2=\sszwei\sim$\hfil}}
\def\Utill_#1{\setbox5=\hbox{$U_#1$}\copy5\kern-\wd5%
\lower4pt\hbox to \wd5 {$\textfont2=\sszwei\sim$\hfil}}
\def\Vtill_#1{\setbox5=\hbox{$V_#1$}\copy5\kern-\wd5%
\lower4pt\hbox to \wd5 {$\textfont2=\sszwei\sim$\hfil}}
\def\Wtill_#1{\setbox5=\hbox{$W_#1$}\copy5\kern-\wd5%
\lower4pt\hbox to \wd5 {$\textfont2=\sszwei\sim$\hfil}}
\def\Xtill_#1{\setbox5=\hbox{$X_#1$}\copy5\kern-\wd5%
\lower4pt\hbox to \wd5 {$\textfont2=\sszwei\sim$\hfil}}
\def\Ytill_#1{\setbox5=\hbox{$Y_#1$}\copy5\kern-\wd5%
\lower4pt\hbox to \wd5 {$\textfont2=\sszwei\sim$\hfil}}
\def\Ztill_#1{\setbox5=\hbox{$Z_#1$}\copy5\kern-\wd5%
\lower4pt\hbox to \wd5 {$\textfont2=\sszwei\sim$\hfil}}

%special combinations with a tilde.

\def\SItil{\mathop{\SI}\limits_{\raise9pt\hbox{$\textfont2=\sszwei\sim$}}\!{}}
\def\Gstill_#1{\setbox5=\hbox{$\scriptstyle {G}_#1$}\copy5\kern-\wd5%
\lower4pt\hbox to \wd5 {$\textfont2=\sszwei\sim$\hfil}}
\def\SItill_#1{\setbox5=\hbox{$\SI_#1$}\copy5\kern-\wd5%
\lower4pt\hbox to \wd5 {$\textfont2=\sszwei\sim$\hfil}}
\def\bStill_#1{\setbox5=\hbox{$\bS_#1$}\copy5\kern-\wd5%
\lower4pt\hbox to \wd5 {$\textfont2=\sszwei\sim$\hfil}}
\def\Qhtill_#1{\setbox5=\hbox{${\widehat Q}_#1$}\copy5\kern-\wd5%
\lower4pt\hbox to \wd5 {$\textfont2=\sszwei\sim$\hfil}}

%{\setbox5=\hbox{$S_{\alpha}$}\copy5\kern-\wd5\lower4pt\hbox to \wd5%
%{$\textfont2=\sszwei\sim$\hfil}}

%{\setbox5=\hbox{$S_i$}\copy5\kern-\wd5\lower4pt\hbox to \wd5%
%{$\textfont2=\sszwei\sim$\hfil}}

%{\setbox5=\hbox{$S_i$}\copy5\kern-\wd5\lower4pt\hbox to \wd5%
%{$\textfont2=\sszwei\sim$\hfil}}

%{\setbox5=\hbox{$Q_i$}\copy5\kern-\wd5\lower4pt\hbox to \wd5%
%{$\textfont2=\sszwei\sim$\hfil}}

%{\setbox5=\hbox{$Q_{\mu}$}\copy5\kern-\wd5\lower4pt\hbox to \wd5%
%{$\textfont2=\sszwei\sim$\hfil}}

%{\setbox5=\hbox{$Q_n$}\copy5\kern-\wd5\lower4pt\hbox to \wd5%
%{$\textfont2=\sszwei\sim$\hfil}}

%{\setbox5=\hbox{$R_i$}\copy5\kern-\wd5\lower4pt\hbox to \wd5%
%{$\textfont2=\sszwei\sim$\hfil}}

%{\setbox5=\hbox{$D_{\alpha}$}\copy5\kern-\wd5\lower4pt\hbox to \wd5%
%{$\textfont2=\sszwei\sim$\hfil}}

%{\setbox5=\hbox{$E_{\alpha}$}\copy5\kern-\wd5\lower4pt\hbox to \wd5%
%{$\textfont2=\sszwei\sim$\hfil}}

%_#1{\setbox5=\hbox{${{\scriptstyle P}_#1$}\copy5\kern-\wd5%
%\lower4pt\hbox to \wd5 {$\textfont2=\sszwei\sim$\hfil}}
\def\bItill_#1{\setbox5=\hbox{$\bI_{#1}$}\copy5\kern-\wd5%
\lower4pt\hbox to \wd5 {$\textfont2=\sszwei\sim$\hfil}}
\def\bStil{\mathop{\bS}\limits_{\raise9pt\hbox{$\textfont2=\sszwei\sim$}}\!{}}
\def\Dptil{\mathop{\Dp}\limits_{\raise9pt\hbox{$\textfont2=\sszwei\sim$}}\!{}}
\def\SAtill_#1{\setbox5=\hbox{$\SA_{#1}$}\copy5\kern-\wd5%
\lower4pt\hbox to \wd5 {$\textfont2=\sszwei\sim$\hfil}}
\def\SAtil{\setbox5=\hbox{$\SA$}\copy5\kern-\wd5%
\lower4pt\hbox to \wd5 {$\textfont2=\sszwei\sim$\hfil}}
\def\cUtill_#1{\setbox5=\hbox{$\SA_{#1}$}\copy5\kern-\wd5%
\lower4pt\hbox to \wd5 {$\textfont2=\sszwei\sim$\hfil}}
\def\cUtil{\mathop{\SA}\limits_{\raise9pt\hbox{$\textfont2=\sszwei\sim$}}\!{}}
\def\cAtill_#1{\setbox5=\hbox{$\SA_{#1}$}\copy5\kern-\wd5%
\lower4pt\hbox to \wd5 {$\textfont2=\sszwei\sim$\hfil}}
\def\cAtil{\mathop{\SA}\limits_{\raise9pt\hbox{$\textfont2=\sszwei\sim$}}\!{}}
\font\funny=cmsy5
\def\supQtil{Q\llap{\lower5pt\hbox{$\textfont2=\funny\sim$\kern0.5pt}}{}}
\def\Dptill_#1{\setbox5=\hbox{$\Dp_#1$}\copy5\kern-\wd5%
\lower4pt\hbox to \wd5 {$\textfont2=\sszwei\sim$\hfil}}

%small letters with a tilde and subset.

\def\atill_#1{\setbox5=\hbox{$a_#1$}\copy5\kern-\wd5%
\lower4pt\hbox to \wd5 {$\textfont2=\sszwei\sim$\hfil}}
\def\btill_#1{\setbox5=\hbox{$b_#1$}\copy5\kern-\wd5%
\lower4pt\hbox to \wd5 {$\textfont2=\sszwei\sim$\hfil}}
\def\ctill_#1{\setbox5=\hbox{$c_#1$}\copy5\kern-\wd5%
\lower4pt\hbox to \wd5 {$\textfont2=\sszwei\sim$\hfil}}
\font\msym=msym10
\def\dtill_#1{\setbox5=\hbox{$d_#1$}\copy5\kern-\wd5%
\lower4pt\hbox to \wd5 {$\textfont2=\sszwei\sim$\hfil}}
\def\etill_#1{\setbox5=\hbox{$e_#1$}\copy5\kern-\wd5%
\lower4pt\hbox to \wd5 {$\textfont2=\sszwei\sim$\hfil}}
\def\ftill_#1{\setbox5=\hbox{$f_#1$}\copy5\kern-\wd5%
\lower6pt\hbox to \wd5 {$\textfont2=\sszwei\sim$\hfil}}
\def\gtill_#1{\setbox5=\hbox{$g_#1$}\copy5\kern-\wd5%
\lower6pt\hbox to \wd5 {$\textfont2=\sszwei\sim$\hfil}}
\def\htill_#1{\setbox5=\hbox{$h_#1$}\copy5\kern-\wd5%
\lower4pt\hbox to \wd5 {$\textfont2=\sszwei\sim$\hfil}}
\def\itill_#1{\setbox5=\hbox{$i_#1$}\copy5\kern-\wd5%
\lower4pt\hbox to \wd5 {$\textfont2=\sszwei\sim$\hfil}}
\def\jtill_#1{\setbox5=\hbox{$j_#1$}\copy5\kern-\wd5%
\lower4pt\hbox to \wd5 {$\textfont2=\sszwei\sim$\hfil}}
\def\ktill_#1{\setbox5=\hbox{$k_#1$}\copy5\kern-\wd5%
\lower4pt\hbox to \wd5 {$\textfont2=\sszwei\sim$\hfil}}
\def\ltill_#1{\setbox5=\hbox{$l_#1$}\copy5\kern-\wd5%
\lower4pt\hbox to \wd5 {$\textfont2=\sszwei\sim$\hfil}}
\def\mtill_#1{\setbox5=\hbox{$m_#1$}\copy5\kern-\wd5%
\lower4pt\hbox to \wd5 {$\textfont2=\sszwei\sim$\hfil}}
\def\ntill_#1{\setbox5=\hbox{$n_#1$}\copy5\kern-\wd5%
\lower4pt\hbox to \wd5 {$\textfont2=\sszwei\sim$\hfil}}
\def\otill_#1{\setbox5=\hbox{$o_#1$}\copy5\kern-\wd5%
\lower4pt\hbox to \wd5 {$\textfont2=\sszwei\sim$\hfil}}
\def\ptill_#1{\setbox5=\hbox{$p_#1$}\copy5\kern-\wd5%
\lower4pt\hbox to \wd5 {$\textfont2=\sszwei\sim$\hfil}}
\def\qtill_#1{\setbox5=\hbox{$q_#1$}\copy5\kern-\wd5%
\lower4pt\hbox to \wd5 {$\textfont2=\sszwei\sim$\hfil}}
\def\rtill_#1{\setbox5=\hbox{$r_#1$}\copy5\kern-\wd5%
\lower4pt\hbox to \wd5 {$\textfont2=\sszwei\sim$\hfil}}
\def\still_#1{\setbox5=\hbox{$s_#1$}\copy5\kern-\wd5%
\lower4pt\hbox to \wd5 {$\textfont2=\sszwei\sim$\hfil}}
\def\ttill_#1{\setbox5=\hbox{$t_#1$}\copy5\kern-\wd5%
\lower4pt\hbox to \wd5 {$\textfont2=\sszwei\sim$\hfil}}
\def\utill_#1{\setbox5=\hbox{$u_#1$}\copy5\kern-\wd5%
\lower4pt\hbox to \wd5 {$\textfont2=\sszwei\sim$\hfil}}
\def\vtill_#1{\setbox5=\hbox{$v_#1$}\copy5\kern-\wd5%
\lower4pt\hbox to \wd5 {$\textfont2=\sszwei\sim$\hfil}}
\def\wtill_#1{\setbox5=\hbox{$w_#1$}\copy5\kern-\wd5%
\lower4pt\hbox to \wd5 {$\textfont2=\sszwei\sim$\hfil}}
\def\xtill_#1{\setbox5=\hbox{$x_#1$}\copy5\kern-\wd5%
\lower4pt\hbox to \wd5 {$\textfont2=\sszwei\sim$\hfil}}
\def\ytill_#1{\setbox5=\hbox{$y_#1$}\copy5\kern-\wd5%
\lower4pt\hbox to \wd5 {$\textfont2=\sszwei\sim$\hfil}}
\def\ztill_#1{\setbox5=\hbox{$z_#1$}\copy5\kern-\wd5%
\lower4pt\hbox to \wd5 {$\textfont2=\sszwei\sim$\hfil}}

%Greek letters with a tilde and subset.

\def\kappatill_#1{\setbox5=\hbox{$\kappa_#1$}\copy5\kern-\wd5%
\lower4pt\hbox to \wd5 {$\textfont2=\sszwei\sim$\hfil}}
\def\sigmatill_#1{\setbox5=\hbox{$\sigma_#1$}\copy5\kern-\wd5%
\lower4pt\hbox to \wd5 {$\textfont2=\sszwei\sim$\hfil}}
\def\lambdatill_#1{\setbox5=\hbox{$\lambda_#1$}\copy5\kern-\wd5%
\lower4pt\hbox to \wd5 {$\textfont2=\sszwei\sim$\hfil}}
\def\mutill_#1{\setbox5=\hbox{$\mu_#1$}\copy5\kern-\wd5%
\lower4pt\hbox to \wd5 {$\textfont2=\sszwei\sim$\hfil}}
\def\thetatill_#1{\setbox5=\hbox{$\theta_#1$}\copy5\kern-\wd5%
\lower4pt\hbox to \wd5 {$\textfont2=\sszwei\sim$\hfil}}
\def\alphatill_#1{\setbox5=\hbox{$\alpha_#1$}\copy5\kern-\wd5%
\lower4pt\hbox to \wd5 {$\textfont2=\sszwei\sim$\hfil}}
\def\gammatill_#1{\setbox5=\hbox{$\gamma_#1$}\copy5\kern-\wd5%
\lower4pt\hbox to \wd5 {$\textfont2=\sszwei\sim$\hfil}}
\def\etatill_#1{\setbox5=\hbox{$\eta_#1$}\copy5\kern-\wd5%
\lower4pt\hbox to \wd5 {$\textfont2=\sszwei\sim$\hfil}}
\def\barxitill_#1{\setbox5=\hbox{$\bar\xi_#1$}\copy5\kern-\wd5%
\lower4pt\hbox to \wd5 {$\textfont2=\sszwei\sim$\hfil}}
\def\rtill_#1{\setbox5=\hbox{$r_#1$}\copy5\kern-\wd5%
\lower4pt\hbox to \wd5 {$\textfont2=\sszwei\sim$\hfil}}
\def\rhotill_#1{\setbox5=\hbox{$\rho_#1$}\copy5\kern-\wd5%
\lower4pt\hbox to \wd5 {$\textfont2=\sszwei\sim$\hfil}}
\def\tautill_#1{\setbox5=\hbox{$\tau_#1$}\copy5\kern-\wd5%
\lower4pt\hbox to \wd5 {$\textfont2=\sszwei\sim$\hfil}}
\def\nutill_#1{\setbox5=\hbox{$\nu_{#1}$}\copy5\kern-\wd5%
\lower4pt\hbox to \wd5 {$\textfont2=\sszwei\sim$\hfil}}
\def\zetatill_#1{\setbox5=\hbox{$\zeta_{#1}$}\copy5\kern-\wd5%
\lower4pt\hbox to \wd5 {$\textfont2=\sszwei\sim$\hfil}}
\def\betatill_#1{\setbox5=\hbox{$\beta_{#1}$}\copy5\kern-\wd5%
\lower4pt\hbox to \wd5 {$\textfont2=\sszwei\sim$\hfil}}
\def\xitill_#1{\setbox5=\hbox{$\xi_{#1}$}\copy5\kern-\wd5%
\lower4pt\hbox to \wd5 {$\textfont2=\sszwei\sim$\hfil}}

% Greek letters with a tilde and a prime.

\def\betaptill_#1{\setbox5=\hbox{$\beta'_{#1}$}\copy5\kern-\wd5%
\lower4pt\hbox to \wd5 {$\textfont2=\sszwei\sim$\hfil}}
\def\tauptill_#1{\setbox5=\hbox{$\tau'_{#1}$}\copy5\kern-\wd5%
\lower4pt\hbox to \wd5 {$\textfont2=\sszwei\sim$\hfil}}
\def\Qptill_#1{\setbox5=\hbox{$Q'_{#1}$}\copy5\kern-\wd5%
\lower4pt\hbox to \wd5 {$\textfont2=\sszwei\sim$\hfil}}
\def\Pptill_#1{\setbox5=\hbox{$\scriptstyle{P}'_{#1}$}\copy5\kern-\wd5%
\lower4pt\hbox to \wd5 {$\textfont2=\sszwei\sim$\hfil}}

% letters with a tilde, subset and superset.

\def\Atills^#1_#2{\setbox5=\hbox{$A^#1_#2$}\copy5\kern-\wd5%
\lower4pt\hbox to \wd5 {$\textfont2=\sszwei\sim$\hfil}}
\def\ttills^#1_#2{\setbox5=\hbox{$t^#1_#2$}\copy5\kern-\wd5%
\lower4pt\hbox to \wd5 {$\textfont2=\sszwei\sim$\hfil}}
\def\alphatills^#1_#2{\setbox5=\hbox{$\alpha^#1_#2$}\copy5\kern-\wd5%
\lower4pt\hbox to \wd5 {$\textfont2=\sszwei\sim$\hfil}}
\def\atills^#1_#2{\setbox5=\hbox{$a^#1_#2$}\copy5\kern-\wd5%
\lower4pt\hbox to \wd5 {$\textfont2=\sszwei\sim$\hfil}}
\def\Ctills^#1_#2{\setbox5=\hbox{$C^#1_#2$}\copy5\kern-\wd5%
\lower4pt\hbox to \wd5 {$\textfont2=\sszwei\sim$\hfil}}
\def\Stills^#1_#2{\setbox5=\hbox{$S^#1_#2$}\copy5\kern-\wd5%
\lower4pt\hbox to \wd5 {$\textfont2=\sszwei\sim$\hfil}}
\def\ptills^#1_#2{\setbox5=\hbox{$p^#1_#2$}\copy5\kern-\wd5%
\lower4pt\hbox to \wd5 {$\textfont2=\sszwei\sim$\hfil}}
\def\zetatills^#1_#2{\setbox5=\hbox{$\zeta^#1_#2$}\copy5\kern-\wd5%
\lower4pt\hbox to \wd5 {$\textfont2=\sszwei\sim$\hfil}}
\def\mutills^#1_#2{\setbox5=\hbox{$\mu^#1_#2$}\copy5\kern-\wd5%
\lower4pt\hbox to \wd5 {$\textfont2=\sszwei\sim$\hfil}}
\def\gtills^#1_#2{\setbox5=\hbox{$g^#1_#2$}\copy5\kern-\wd5%
\lower6pt\hbox to \wd5 {$\textfont2=\sszwei\sim$\hfil}}
\def\ftills^#1_#2{\setbox5=\hbox{$f^#1_#2$}\copy5\kern-\wd5%
\lower6pt\hbox to \wd5 {$\textfont2=\sszwei\sim$\hfil}}
\def\ntills^#1_#2{\setbox5=\hbox{$n^#1_#2$}\copy5\kern-\wd5%
\lower4pt\hbox to \wd5 {$\textfont2=\sszwei\sim$\hfil}}
\def\mtills^#1_#2{\setbox5=\hbox{$m^#1_#2$}\copy5\kern-\wd5%
\lower4pt\hbox to \wd5 {$\textfont2=\sszwei\sim$\hfil}}
\def\Qtills^#1_#2{\setbox5=\hbox{$Q^#1_#2$}\copy5\kern-\wd5%
\lower4pt\hbox to \wd5 {$\textfont2=\sszwei\sim$\hfil}}
\def\barQtills^#1_#2{\setbox5=\hbox{$\bar{Q}^#1_#2$}\copy5\kern-\wd5%
\lower4pt\hbox to \wd5 {$\textfont2=\sszwei\sim$\hfil}}
\def\Ftills^#1_#2{\setbox5=\hbox{$F^#1_#2$}\copy5\kern-\wd5%
\lower4pt\hbox to \wd5 {$\textfont2=\sszwei\sim$\hfil}}
\def\xitills^#1_#2{\setbox5=\hbox{$\xi^#1_#2$}\copy5\kern-\wd5%
\lower4pt\hbox to \wd5 {$\textfont2=\sszwei\sim$\hfil}}
\def\Gtills^#1_#2{\setbox5=\hbox{$G^#1_#2$}\copy5\kern-\wd5%
\lower4pt\hbox to \wd5 {$\textfont2=\sszwei\sim$\hfil}}
\def\qtills^#1_#2{\setbox5=\hbox{$q^#1_#2$}\copy5\kern-\wd5%
\lower4pt\hbox to \wd5 {$\textfont2=\sszwei\sim$\hfil}}
\def\Ttills^#1_#2{\setbox5=\hbox{$T^#1_#2$}\copy5\kern-\wd5%
\lower4pt\hbox to \wd5 {$\textfont2=\sszwei\sim$\hfil}}
\def\betatills^#1_#2{\setbox5=\hbox{$\beta^#1_#2$}\copy5\kern-\wd5%
\lower4pt\hbox to \wd5 {$\textfont2=\sszwei\sim$\hfil}}
\def\cItills^#1_#2{\setbox5=\hbox{$\dbI^#1_#2$}\copy5\kern-\wd5%
\lower4pt\hbox to \wd5 {$\textfont2=\sszwei\sim$\hfil}}

% special letters, each with different names throughout the file.

\def\BBD{\shD}

% special signs and operations

\def\notVdash{\rlap/\kern-2pt\Vdash}
\def\notvdash{\rlap/\kern-6pt\vdash}
\def\eqdf{\buildrel\rm def\over =}
\def\notVdash{\rlap/\kern-2pt\Vdash}
\def\lc{\mathrel{<\mkern-11mu\raise1.05pt\hbox{$\scriptstyle\circ$}}}

\def\lesseqdot{\mathrel{\le\mkern-11mu\raise1.85pt\hbox{$\scriptstyle\circ$}}}

\def\lapprox{\rlap/\kern-0pt\approx}
\def\smallbox#1{\leavevmode\thinspace\hbox{\vrule\vtop{\vbox
   {\hrule\kern1pt\hbox{\vphantom{\tt/}\thinspace{\tt#1}\thinspace}}
   \kern1pt\hrule}\vrule}\thinspace}

% roman words

\def\Dp{{\rm Dp}}
\def\lg{\,{\ell g}}

\def\Min{\hbox{\norm Min}}

\def\Min{\hbox{\rm Min}}

\def\max{{\rm \ max\,}}
\def\mod{{\rm \ mod\ }}

\def\tp{{\rm \ tp}}

\def\min{{\rm \ min \,}}

\def\Rang{{\rm Rang}}

\def\log {\rm log}

%\def\underW{\bf W}

% small letters with a tilde and no sub/superset.

\def\atil{\mathop{a}\limits_{\raise9pt\hbox{$\textfont2=\sszwei\sim$}}\!{}}
\def\btil{\mathop{b}\limits_{\raise9pt\hbox{$\textfont2=\sszwei\sim$}}\!{}}
\def\ctil{\mathop{c}\limits_{\raise9pt\hbox{$\textfont2=\sszwei\sim$}}\!{}}
\def\dtil{\mathop{d}\limits_{\raise9pt\hbox{$\textfont2=\sszwei\sim$}}\!{}}
\def\etil{\mathop{e}\limits_{\raise9pt\hbox{$\textfont2=\sszwei\sim$}}\!{}}
\def\ftil{\mathop{f}\limits_{\raise9pt\hbox{$\textfont2=\sszwei\sim$}}\!{}}
\def\gtil{\mathop{g}\limits_{\raise9pt\hbox{$\textfont2=\sszwei\sim$}}\!{}}
\def\htil{\mathop{h}\limits_{\raise9pt\hbox{$\textfont2=\sszwei\sim$}}\!{}}
\def\itil{\mathop{i}\limits_{\raise9pt\hbox{$\textfont2=\sszwei\sim$}}\!{}}
\def\jtil{\mathop{j}\limits_{\raise9pt\hbox{$\textfont2=\sszwei\sim$}}\!{}}
\def\ktil{\mathop{k}\limits_{\raise9pt\hbox{$\textfont2=\sszwei\sim$}}\!{}}
\def\ltil{\mathop{l}\limits_{\raise9pt\hbox{$\textfont2=\sszwei\sim$}}\!{}}
\def\mtil{\mathop{m}\limits_{\raise9pt\hbox{$\textfont2=\sszwei\sim$}}\!{}}
\def\ntil{\mathop{n}\limits_{\raise9pt\hbox{$\textfont2=\sszwei\sim$}}\!{}}
\def\otil{\mathop{o}\limits_{\raise9pt\hbox{$\textfont2=\sszwei\sim$}}\!{}}
\def\ptil{\mathop{p}\limits_{\raise9pt\hbox{$\textfont2=\sszwei\sim$}}\!{}}
\def\qtil{\mathop{q}\limits_{\raise9pt\hbox{$\textfont2=\sszwei\sim$}}\!{}}
\def\rtil{\mathop{r}\limits_{\raise9pt\hbox{$\textfont2=\sszwei\sim$}}\!{}}
\def\stil{\mathop{s}\limits_{\raise9pt\hbox{$\textfont2=\sszwei\sim$}}\!{}}
\def\ttil{\mathop{t}\limits_{\raise9pt\hbox{$\textfont2=\sszwei\sim$}}\!{}}
\def\vtil{\mathop{v}\limits_{\raise9pt\hbox{$\textfont2=\sszwei\sim$}}\!{}}
\def\wtil{\mathop{w}\limits_{\raise9pt\hbox{$\textfont2=\sszwei\sim$}}\!{}}
\def\xtil{\mathop{x}\limits_{\raise9pt\hbox{$\textfont2=\sszwei\sim$}}\!{}}
\def\ytil{\mathop{y}\limits_{\raise9pt\hbox{$\textfont2=\sszwei\sim$}}\!{}}
\def\ztil{\mathop{z}\limits_{\raise9pt\hbox{$\textfont2=\sszwei\sim$}}\!{}}
\def\bbitil{\mathop{\dbI}\limits_{\raise9pt\hbox{$\textfont2=\sszwei\sim$}}\!{}}
\def\alphaitil{\mathop{\alpha(i)}\limits_{\raise9pt\hbox{$\textfont2=\sszwei\sim$}}\!{}}

% Bold face letters

\font\fatone=cmssbx10
\def\bS{\hbox{\fatone S}}
\def\bI{\hbox{\fatone I}}

\def\bt{\hbox{\fatone t}}
\def\bT{\hbox{\fatone T}}

% capital letters with a tilde and no sub/superset

\def\Atil{A\llap{\lower5pt\hbox{$\textfont2=\sszwei\sim$\kern2pt}}{}}
\def\Btil{B\llap{\lower5pt\hbox{$\textfont2=\sszwei\sim$\kern2pt}}{}}
\def\Ctil{C\llap{\lower5pt\hbox{$\textfont2=\sszwei\sim$\kern2pt}}{}}
\def\Dtil{D\llap{\lower5pt\hbox{$\textfont2=\sszwei\sim$\kern2pt}}{}}
\def\Etil{E\llap{\lower5pt\hbox{$\textfont2=\sszwei\sim$\kern2pt}}{}}
\def\Ftil{F\llap{\lower5pt\hbox{$\textfont2=\sszwei\sim$\kern2pt}}{}}
\def\Gtil{G\llap{\lower5pt\hbox{$\textfont2=\sszwei\sim$\kern2pt}}{}}
\def\Htil{H\llap{\lower5pt\hbox{$\textfont2=\sszwei\sim$\kern2pt}}{}}
\def\Itil{I\llap{\lower5pt\hbox{$\textfont2=\sszwei\sim$\kern2pt}}{}}
\def\Jtil{J\llap{\lower5pt\hbox{$\textfont2=\sszwei\sim$\kern2pt}}{}}
\def\Ktil{K\llap{\lower5pt\hbox{$\textfont2=\sszwei\sim$\kern2pt}}{}}
\def\Ltil{L\llap{\lower5pt\hbox{$\textfont2=\sszwei\sim$\kern2pt}}{}}
\def\Mtil{M\llap{\lower5pt\hbox{$\textfont2=\sszwei\sim$\kern2pt}}{}}
\def\Ntil{N\llap{\lower5pt\hbox{$\textfont2=\sszwei\sim$\kern2pt}}{}}
\def\Otil{O\llap{\lower5pt\hbox{$\textfont2=\sszwei\sim$\kern2pt}}{}}
\def\Ptil{P\llap{\lower5pt\hbox{$\textfont2=\sszwei\sim$\kern2pt}}{}}
\def\Qtil{Q\llap{\lower5pt\hbox{$\textfont2=\sszwei\sim$\kern2pt}}{}}
\def\Rtil{R\llap{\lower5pt\hbox{$\textfont2=\sszwei\sim$\kern2pt}}{}}
\def\Stil{S\llap{\lower5pt\hbox{$\textfont2=\sszwei\sim$\kern2pt}}{}}
\def\Ttil{T\llap{\lower5pt\hbox{$\textfont2=\sszwei\sim$\kern2pt}}{}}
\def\Util{U\llap{\lower5pt\hbox{$\textfont2=\sszwei\sim$\kern2pt}}{}}
\def\Vtil{V\llap{\lower5pt\hbox{$\textfont2=\sszwei\sim$\kern2pt}}{}}
\def\Wtil{W\llap{\lower5pt\hbox{$\textfont2=\sszwei\sim$\kern2pt}}{}}
\def\Xtil{X\llap{\lower5pt\hbox{$\textfont2=\sszwei\sim$\kern2pt}}{}}
\def\Ytil{Y\llap{\lower5pt\hbox{$\textfont2=\sszwei\sim$\kern2pt}}{}}
\def\Ztil{Z\llap{\lower5pt\hbox{$\textfont2=\sszwei\sim$\kern2pt}}{}}
\def\Qtil{Q\llap{\lower5pt\hbox{$\textfont2=\sszwei\sim$\kern2pt}}{}}
\def\Qtil{Q\llap{\lower5pt\hbox{$\textfont2=\sszwei\sim$\kern2pt}}{}}
\def\Qtil{Q\llap{\lower5pt\hbox{$\textfont2=\sszwei\sim$\kern2pt}}{}}
\def\Qtil{Q\llap{\lower5pt\hbox{$\textfont2=\sszwei\sim$\kern2pt}}{}}
\def\Qtil{Q\llap{\lower5pt\hbox{$\textfont2=\sszwei\sim$\kern2pt}}{}}
\def\Stil{\mathop{S}\limits_{\raise9pt\hbox{$\textfont2=\sszwei\sim$}}\!{}}
\def\Vtil{\mathop{V}\limits_{\raise9pt\hbox{$\textfont2=\sszwei\sim$}}\!{}}
\def\Ttil{\mathop{T}\limits_{\raise9pt\hbox{$\textfont2=\sszwei\sim$}}\!{}}
\def\Wtil{\mathop{W}\limits_{\raise9pt\hbox{$\textfont2=\sszwei\sim$}}\!{}}
\def\Xtil{\mathop{X}\limits_{\raise9pt\hbox{$\textfont2=\sszwei\sim$}}\!{}}
\def\Ytil{\mathop{Y}\limits_{\raise9pt\hbox{$\textfont2=\sszwei\sim$}}\!{}}
\def\Ztil{\mathop{Z}\limits_{\raise9pt\hbox{$\textfont2=\sszwei\sim$}}\!{}}

% greek and special letters with tilde and no sub/superset.

\def\mutil{\mathop{\mu}\limits_{\raise9pt\hbox{$\textfont2=\sszwei\sim$}}\!{}}
\def\gammatil{\mathop{\gamma}\limits_{\raise9pt\hbox{$\textfont2=\sszwei\sim$}}\!{}}
\def\tautil{\mathop{\tau}\limits_{\raise9pt\hbox{$\textfont2=\sszwei\sim$}}\!{}}
\def\zetalil{\mathop{\zeta}\limits_{\raise9pt\hbox{$\scriptscriptstyle\textfont2=\sszwei\sim$}}\!{}}
\def\rhotil{\mathop{\rho}\limits_{\raise9pt\hbox{$\textfont2=\sszwei\sim$}}\!{}}
\def\nutil{\mathop{\nu}\limits_{\raise9pt\hbox{$\textfont2=\sszwei\sim$}}\!{}}
\def\xitil{\mathop{\xi}\limits_{\raise9pt\hbox{$\textfont2=\sszwei\sim$}}\!{}}
\def\Xitil{\mathop{\Xi}\limits_{\raise9pt\hbox{$\textfont2=\sszwei\sim$}}\!{}}
\def\barxitil{\mathop{\bar\xi}\limits_{\raise9pt\hbox{$\textfont2=\sszwei\sim$}}\!{}}
\def\etatil{\mathop{\eta}\limits_{\raise9pt\hbox{$\textfont2=\sszwei\sim$}}\!{}}
\def\BBDtil{\mathop{\BBD}\limits_{\raise9pt\hbox{$\textfont2=\sszwei\sim$}}\!{}}
\def\zetatil{\mathop{\zeta}\limits_{\raise9pt\hbox{$\textfont2=\sszwei\sim$}}\!{}}
\def\thetatil{\mathop{\theta}\limits_{\raise9pt\hbox{$\textfont2=\sszwei\sim$}}\!{}}

% ``Double'' Letters: Like the signs for rational/natural numbers.

\font\msym=msym10
\textfont9=\msym
\font\msyms=msym7
\scriptfont9=\msyms

\def\dbI{{\fam9 I}}
\def\dbR{{\fam9 R}}

\def\dbD{{\fam9 D}}
\def\BBD{\dbD}

\def\BBR{\dbR}

\def\bbitills^#1_#2{\setbox5=\hbox{$\dbI^#1_#2$}\copy5\kern-\wd5%
\lower4pt\hbox to \wd5 {$\textfont2=\sszwei\sim$\hfil}}
\def\cItill_#1{\setbox5=\hbox{$\dbI_#1$}\copy5\kern-\wd5%
\lower4pt\hbox to \wd5 {$\textfont2=\sszwei\sim$\hfil}}
\def\bbitill_#1{\setbox5=\hbox{$\dbI_#1$}\copy5\kern-\wd5%
\lower4pt\hbox to \wd5 {$\textfont2=\sszwei\sim$\hfil}}

% Gothic Letters

\font\eufb=eufb10
\def\GA{\hbox{\eufb\char"41}}
\def\GB{\hbox{\eufb\char"42}}

\def\Gt{\hbox{\eufb\char"74}}
\font\cmbx=cmbx12
\font\scmbx=cmbx7
\def\smalltriangle{\hbox{\scmbx\char"1}}
\def\bigtriangle{\hbox{\cmbx\char"1}}

\hfuzz72pt
\def\rkol{0-1 Law}%, \today}
\def\lkol{0-1 Law}%, \today}
\def\|{\Big|}
\def\fo{{\rm fo}}

\def\LL{\SL}
\def\opo{{\rm opo}}
\def\<{\langle}
\def\>{\rangle}
\def\restrictedto{\uhr}

\font\eufb=eufb10
\def\GA{\hbox{\eufb\char"41}}
\def\GB{\hbox{\eufb\char"42}}

\def\bbn{{\fam9 N}}
\def\bbN{{\bbn}}
\def\bbz{{\fam9 Z}}
\def\bbr{{\fam9 R}}

\def\vare{\varepsilon}
\def\yes{{\rm yes}}
\def\if{{\rm if\,}}

\def\og{{\rm og}}
\def\no{{\rm no}}

\def\UTH{{\rm UTH}}

\def\or{{\rm or}}

\def\npr{{\rm npr}}
\def\Rang{{\rm Rang\,}}
\def\th{{\rm th}}
\def\bth{{\rm bth}}
\def\TH{{\rm TH}}
\def\BTH{{\rm BTH}}
\def\uth{{\rm uth}}

\def\nl{\nz}
\def\spr{{\rm spr}}

\def\id{{\rm id}}
\def\Prob{{\rm Prob}}
\def\dis{{\rm dis}}

\def\tp{{\rm tp}}
\def\qf{{\rm qf}}

\def\sut{{\rm sut}}
\def\rsim{{\rm sim}}
\def\ws{{\rm ws}}
\def\as{{\rm as}}

\def\title{On the Very Weak 0-1 Law for Random Graphs with Orders}
\def\abstract{ Let us draw a graph $R$ on
$\{0,1,\dots,n-1\}$ by having an edge $\{i,j\}$ with probability
$p_{|i-j|}$, where $\sum_i p_i<\infty$, and let $M_n=(n, <, R)$.  For a
first order sentence $\psi$ let $a^n_\psi$ be the probability of $M_n\models
\psi$.  We know that the sequence $a^1_\psi, a^2_\psi, \dots, a^n_\psi,
\dots$ does not necessarily converge. But here we find a weaker substitute
which we call the very weak 0-1 law.  We prove that $\lim_{n\to \infty}
(a^n_\psi-a^{n+1}_\psi)=0$. For this we need a theorem on the (first order)
theory of distorted sum of models.} 

\def\thanks {{ \small\par \openup-5pt
Research partially supported by the Binational Science Foundation and
partially sponsored by the Edmund Landau Center for research in Mathematical
Analysis, supported by the Minerva Foundation (Germany).  Done fall 91,
Publication 463.\par}}
\def\shelahfirst{{\topskip=1.5 true in \centerline\shelahtitle
\vskip 1 true in \itemitem{{\sl Abstract\/}:}
\abstract\par\vfill\noindent\thanks\par}\break}
\def\shelahtitle{\vbox{\hbox{\ta{}{\title}}\vskip.6true in
\hbox{Saharon Shelah}\medskip\hbox{\sl The Hebrew University, Math Institute}
\hbox{\sl Rutgers University, Math Department}}}
\shelahfirst

%%% \input 463s0
%
% file 463s0.tex starts here 
%
\tb{\S0}{Introduction}

The kind of random models $M_n=(n, <,R)$ from the abstract are from Luczak
Shelah [LuSh 435] where among other things, it is proved that the
probability $a^n_\psi=:\Prob(M_n\models \psi$) of $M_n\models \psi$ does not
necessarily converge (but if $\sum_i ip_i<\infty$ then it converges, and the
value of $\sum p_i$ is not enough, in general, to determine convergence).
The theorem in the abstract appears in \S1 and is proved in \S3, it
says that the sequence of probabilities still behaves
(somewhat) nicely.

The first results (in various probabilistic
distributions) on the asymptomatic behavior of $a^n_\psi$
(see Glebski et al [GKLT], Fagin [F] and survey [Sp]) say that it
converges to $0$ or $1$ hence the name zero one law.  In other cases weaker
results were gotten: $a^n_\psi$ converges (to some real). We suggest an even
weaker version: $\vert a^{n+1}_\psi-a^n_\psi\vert $ converges to zero.
We also define $h$-very weak zero one law (see Definition 1.2(2)), but
concentrate on the one above.
Note
that many examples to nonconvergence are done by finding $\psi$ such that
e.g. if $\log(n)=1\mod 10$ then $a^n_\psi\sim 1$ and if $\log(n)=5\mod 10$
then $a^n_\psi\sim 5$ (or even using functions with $h(n)\to \infty$,
$h(n)<<\log(n)$).  As the most known results were called zero one law we
prefer the name ``very weak 0-1 law'' on weak convergence.  A (first order)
sentence whose probability $a^n_\psi$ (defined above, for the distribution
defined above) may not converge (by [LuSh 435]) is $\psi_0=:(\exists x)
\forall yz(y<x\le z\to \neg yRz).$ But we can find a sequence
$m_0<m_1<\dots$ such that the probability that $\Phi=:\bigvee_i (\exists
yz)(y\le m_i\, \&\, m_{i+1}\le z\, \&\, yRz)$ is very small.  How do we
prove $0=\lim (a^n_\psi-a^{n+1}_\psi)$? By changing the rule of making the
random choice we get $M'_n$ nice enough, ensuring $\Phi$ holds, while the
probability changes little (see \S3). Now $M'_n$ is almost the sum of
$M'_n\uhr (m_i, m_{i+1})$, precisely $M'_n$ is determined by $M'_n\uhr[m_i,
m_{i+2})$, for $i=0, 1, 2,\dots$, so we call it a distorted sum.  Now a
model theoretic lemma from \S2 on the $n$-theory of a distorted sum of
models enables us to prove the main theorem in \S3 (on the model
theoretic background see \S2). Later in \S4, \S5 we
deal with some refinements not needed for the main theorem (1.4).

In \S1 we also get the very weak 0-1 law for a random partial order
suggested by Luczak.  In a subsequent paper [Sh 548] we prove the very
weak zero law
for some other very natural cases: e.g. for a random 2-place function and
for $(n, <, R)$ with $<$ the natural order, $R$ a random graph (=symmetric
irreflexive relation) with edge probability $p$. In another one, [Sh
467] we deal with zero one law for the random model from the abstract
with $p_i={1\over{i^a}}$ (mainly: no order, $a\in (0,1)$ irrational).
See also [Sh 550], [Sh 551], [Sh 581]. Spencer is continuing [Sh 548]
looking at the exact $h$ for which $h$-very weak zero one law holds
(see Definition 1.2(2) here).
I thank Shmuel
Lifsches for many corrections.

\lz {\bf Notation} \nl
\item{\ } $\bbn$ is the set of natural numbers.
\item{\ } We identify $n\in \bbn$ with the set $\{0, \ldots, n-1\}$.
\item{\ } $\bbz$ is the set of integers.
\item{\ } $\bbr$ is the set of reals, $\bbr^+$ is the set of reals
which are positive (i.e. $>0$).
\item{\ } $i,j,k, \ell, m,n,r,s,t$ are natural numbers.
\item{\ } $\vare, \zeta$ are positive reals (or functions
with values in $\bbr^+$).
\item{\ } $f, g, h$ are functions.
\item{\ } $\tau$ denotes a vocabulary (for simplicity- set of predicates),
$\SL^{\fo}_\tau$ is the set of first
 order formulas in the vocabulary $\tau$.
(Generally, if $\SL$ is a logic $\SL_\tau$ the set of
 sentences (or formulas) in $\SL$ in the vocabulary $\tau$
and $\SL^{{\rm fo}}$ is the first order logic).
\item{\ } For a first order sentence (or formula) $\varphi$ let $d_\varphi=d[\varphi]$
be its quantifier depth.
\item{\ } $M,N$ denote models, but we do not distinguish
strictly between a model and its universe = set of elements.
\item{\ } $\tau(M)$ is the vocabulary of $M$, for $R\in\tau(M)$,
$n(R)$ is the number of places (=arity of $R$),
$R^M$ the interpretation of $R$ in $M$.
\item{\ } A basic formula is one of the form
$\pm R(x_{i_0}, \dots, x_{i_{n(R)-1}})$
(i.e. $R(x_{i_0}, \dots, x_{i_{n(R)-1}})$ or $\neg R(x_{i_0}, \dots,
x_{i_{n(R)-1}})$
\item{\ } $\bar a$ denotes a sequence of elements of a model.
$\lg(\bar a)$ is the length of $\bar a$.

If $<$ belongs to the vocabulary $\tau$, then in $\tau$-models $<^M$ is a
linear order, if not said otherwise.  If $M_i$ are $\tau$-models for $i<n$,
$M=\sum_{i<n}M_i$ is (assuming for simplicity the universes are pairwise
disjoint) the models defined by: universe $\bigcup_{i<n}M_i,
R^M=\bigcup_{i<n} R^{M_i}$ for $R\in \tau$ except that if $<\;\in \tau$
then: $x<^M y\Leftrightarrow \bigvee_{i<j<n}[x\in M_i\,\&\, y\in M_j]\,\vee\,
\bigvee_{i<n} x<^{M_i}y$ (similarly with any linear order $I$ instead of
$n$).  We write $M_0+M_1$ instead of $\sum_{i<2} M_i$.\nl
$\psi^{{\rm if(\theta)}}$ is $\psi$ if $\theta$ is true, $\neg \psi$ if
$\theta$ is false.
\nl
We identify true, false with yes, no.
\nl
Note: $t$ is a natural number,  $\bt$ kind of depth $n$ theory, $\Gt$ is a
truth value, $\smalltriangle$ is the symmetric difference, $\Delta$
denote a set of formulas.

%
% file 463s0.tex ends  here 
%
%%% \input 463s1
%
% file 463s1.tex starts here 
%
\tb{\S1}{The Very Weak Zero One Law.}

\noindent {\bf 1.1 Definition:} 1) A 0-1 law context is a sequence $\bar
K=\langle K_n, \mu_n:n<\omega\rangle$ such that: \itemitem{(a)} for some
vocabulary $\tau=\tau_{\bar K}$, for every $n$, $K_n$ is a family of
$\tau$-models, closed under isomorphism with the family of isomorphism types
being a set.
\itemitem{(b)} for each $n$, $\mu_n$ is a probability measure on the set of
isomorphism types (of models from $K_n$).  \item{2)} $\bar K$ is finitary
(countable) if for each $n$ the set $\{M/\cong\;: M\in K_n\}$ is finite
(countable).
\item{3)} For a sentence $\psi$ (not necessarily f.o.)
$\Prob_{\mu_n}(M_n\models\psi)$ or $\Prob_{K_n}(M_n\models \psi)$ or
$\Prob_{\mu_n}(M_n\models \psi \| M_n\in K_n$) means $\mu_n\{M/{\cong}\,:
M\in K_n, M\models\psi\}$.
\item{4)} Instead of clause (a) of (1), we may use $K_n$ a set of
$\tau$-models, $\mu_n$ a probability measure on $K_n$; particularly we
introduce a random choice of $M_n$; the translation between the two contexts
should be clear.
\lz
{\bf 1.1A Discussion:} A 0-1 law context \underbar{is not} necessarily a
context which satisfies a 0-1 law. It is a context in which we can
formulate a 0-1 law, and also weaker variants.
\lz
{\bf 1.2 Definition:}
0) A 0-1 context $\bar K$ satisfies the 0-1 law for a logic $\SL$
\underbar{if} for every sentence $\varphi\in \SL_\tau$ (with
$\tau=\tau_{\bar K}$ of course) we have:
$a^n_\varphi\eqdf \Prob_{\mu_n}(M_n\models \varphi)$ converges to zero or
converges to 1 when $n\to \infty$.
\nz
1) $\bar K$ satisfies the very weak 0--1 law for the logic $\SL$
\underbar{if } for every sentence $\varphi\in \SL_\tau$ we have:
$a_n\eqdf\Prob_{\mu_{n}}(M_{n+1}\models \varphi \| M_{n+1}\in K_{n+1})-
\Prob_{\mu_{n}}(M_{n}\models \varphi \|M_{n}\in K_{n})$ converges to zero as
$n\to \infty$.  \nl
2) $\bar K$ satisfies the $h$-very 0--1 one law for the logic $\SL$
\underbar{if} for every sentence $\varphi\in \SL_\tau$,
$\max_{m_1, m_2\in[n, h(n)]} \Big[ \Prob_{\mu_{m_1}}(M_{m_1}\models \varphi
\| M_{m_1}\in K_{m_1})- \Prob_{\mu_{m_2}}(M_{m_2}\models \varphi\|M_{m_2}\in
K_{m_2})\Big]$ converge to zero as $n\to \infty$.  (We shall concentrate on
part 1).  \nl
3) $\bar K$ satisfies the convergence law for $\SL$ if for every $\varphi\in
\SL_\tau$ we have: $\lk \Prob_{\mu_n} ( M_n\models\varphi\|M_n\in
K_n):n<\omega\rangle$ converges (to some real $\in[0,1]_{\bbr}$).  (So if it
always converges to $0$ or to $1$ we say that $\bar K$ satisfies the 0-1 law
for $\SL$). \nl
4) If $\SL$ is the first order logic we may omit it.

The following $\bar K$ is from Luczak Shelah [LuSh 435].  \lz

{\bf 1.3 Definition:} Let $\bar p=\langle p_i:i\in \bbn\rangle$, $p_i$ a
real number $0\le p_i\le 1$, $p_0=0$.  We define $\bar K^{\og}_{\bar p} $ as
follows: \nl
the models from $K_n$ are of the form $M=(n, <, R)$, so we are using the
variant from\footnote{${}^\dagger$}{ but no two models are isomorphic so the
difference is even more trivial} Definition 1.1(4) where $n=\{0, \dots,
n-1\}$, $<$ the usual order, $R$ a graph (i.e. a symmetric relation on
$n$ which is irreflexive i.e.
$\neg xRx$) and: $\Prob_{\mu_n}(M_n\cong M\| M_n\in K_n)$ is
$\prod\{p_{j-i}:i<j<n$ and $iRj\}$ $\times\prod\{1- p_{j-i}:i<j<n$ and $\neg
iRj\}$, i.e. for each $i<j<n$ we decide whether $iRj$ by flipping a coin
with probability $p_{j-i} $ for yes, independently for the different pairs.
\lz

{\bf 1.4 Theorem:} $\bar K^{\og}$ satisfies the very weak 0--1 law
\underbar{if} $\sum_i p_i<\infty$. \nz

This is our main result. We shall prove it later in \S3. Luczak [Lu]
suggested another context: \lz

{\bf 1.5 Definition:} $\bar K^{\opo}_{p_n}$ is the following 0-1 law context
$(p_n$ is a function of $n$, $0\le p_n\le 1$; so more precisely we are
defining $\lk K^\opo_{p_n}:n<\omega\rk$).  The models are of the form, $(n,
<, <^*)$, $n=\{0, 1, \dots, n-1\}$, $<$-the usual order, $<^*$-the following
partial order: we draw a random graph on $n$ (edge relation $\bbr$) with
edge probability $p_n$, $x<^*y$ iff there are $k\in N$ and
$x=x_0<x_1<\dots<x_k=y<n$ such that $x_\ell Rx_{\ell+1}$.  The probability is
derived from the probability distribution for $R$ (which does not appear in
the model.)  \lz

{\bf 1.6 Theorem} Assume $p_n={1\over (n+1)^a}$, where $0<a<1$ (like [ShSp
304]).  Then $\bar K^{opo}_{p_n}$ satisfies the very weak zero one law.  \hz

{\it Proof:} Similar to the previous theorem as easily \nl (*) for $\vare>0$
we have: $\Prob(M_n\models(\exists x<y) [\neg x<^* y \& x+ n^{1-\vare}\le y
$) is very small.
\hfill\qed$_{1.6}$

%
% file 463s1.tex ends  here 
%
%%% \input 463s2
%
% file 463s2.tex starts here 
%
\tb {\S2}{ Model Theory: Distorted sum of Models}

The main lemma (2.14) generalizes the addition theory and deals with
models with 
distances (both from model theory).  Concerning the first, see Feferman
Vaught [FV]. The method has its origin in Mostowski [Mo], who dealt with
reduced products. The first work on other products is Beth [B] who dealt
with the ordinal sum of finitely many ordered systems. For a presentation and
history see Feferman Vaught [FV], pp 57--59 and Gurevich [Gu] (particularly
the $\th^n$'s).  Concerning models with distance see Gaifman [Gf], a
forerunner of which was Marcus [M1] who deals with the case of $M=(M, F,
P_i)_{i< n}$, $F$ a unary function, $P_i$ unary predicates and the distance
is as in the graph which the function $F$ defines (i.e. $x, y$
connected by an edge if:
$x=F(y)$ or $y=F(x)$; where for a graph $G$ the distance is $d_G(x,y)=\min
\{k:$ we can find $x_0, \dots, x_k$ such that: $x=x_0$, $y=x_k$, and
$x_\ell, x_{\ell+1}$ are connected$\}$).  We may look at our subject
here as dealing with
sums with only local ``disturbances'', ``semi sums''; ``distorted sums''.
The connections are explained in \S4, \S5.

In 2.16 we draw the conclusion for linear order which we shall use in
\S3 to prove theorem 1.4, in fact proving the main theorem 1.4 from \S2
we use almost only 2.16. Elsewhere we shall return to improving the
numerical bound involved in the proof see [Sh, F-120]
\nl
{\bf Note:} $\GB$ is a (two sorted) model.
\nl
\lz
{\bf 2.1 Definition:} 1) We call $\sigma$ a vocabulary of systems if
$\sigma=\langle \tau_1, \tau_2\rangle$, $\tau_1, \tau_2$ are sets of
predicates (usually finite but not needed).\nl
2) We call $\GB$ a $\sigma$-system if:
\item{(A)} $\GB=(M, I, h, d)\; (=(M^{\GB}, I^{\GB}, h^{\GB}, d^{\GB}))$.
\item{(B)} $M$ is a $\tau_1$-model and $I$ is a $\tau_2$-model,
but we use $M$, $I$
also for their universes.
\item{(C)} $h$ is a function from $M$ onto $I$,
\item{(D)} $d$ is a distance function on $I$, i.e.
\itemitem{($\alpha$)} $d$ is a symmetric two place function from $I$ to
$\bbn\cup\{\infty\} = \{0, 1, 2, 3, \dots, \infty\}$,
\itemitem{$(\beta)$} $d(x, x)=0$
\itemitem{$(\gamma)$} $d(x,z)\le d(x, y)+d(y, z)$.

\noindent
\item{(E)} $M, I$ are
disjoint.

\noindent
3) Let $\sigma(\GB)=\sigma$ for $\GB$ a $\sigma$-system and
$\tau_\ell(\GB)=\tau_\ell$ if $\sigma=\langle \tau_1, \tau_2\rangle$.
\lz
{\bf 2.1A Discussion:} The demands $\lq M, I$ are disjoint and $h$ is onto
$I\rq$ are not essential, this is just for convenience of presentation.  If
$h$ is not onto $I$, we should allow relations on $M\cup I$, but then if $M,
I$ are not disjoint then for each predicate, for each of its places we
should assign: does there appear a member of $M$ or a member of $I$; also we
should have two kinds of variables not allowing equality between variables
of the different kinds.  So in application we may use those alternative
presentations.
\lz
{\bf 2.2 Conventions:} (1) We may allow function symbols (and individual
constants) but then treat them as relations.\nl
(2) We stipulate $h(x)=x$ for $x\in I$ (remember 2.1(2)(E)).\nl
(3) $\bar a\subseteq \GB$ means $\Rang(\bar a)\subseteq M^{\GB}\cup
I^{\GB}$.\nl 
(4) A model $M$ will be identified with the system $\GB=\GB\, {}^\rsim [M]:
M^{\GB}=M$, $I^{\GB}=M$, $h^{\GB}=\id_M$,
$$d(x, y)=\cases{0& $x=y$\cr\infty& $x\not= y$\cr}$$
(you may of course take two disjoint copies of $M$ as $M^{\GB}$ and $I^{\GB}$).
\nl
(5) From a model $M$ we can derive another system $\GB^\dis [M]$:
$M^{\GB}=M$,$I^{\GB}=M$, $h^{\GB}=\id_M$, and $d(x, y)\eqdf\Min \{n:$ there
are $z_0, \dots, z_n\in M$, $x=z_0, y=z_n$ and for $\ell<n$ for some $R\in
\tau(M)$, and sequence $\bar a\in R^M$ we have
$\{z_\ell,z_{\ell+1}\}\subseteq \Rang \bar a\}$ (remember $\tau(M)$ is the
vocabulary of $M$; this is the definition of distance in Gaifman [Gf]).
\lz
{\bf 2.2A Remark:} 
1) So the difference between $\GB^{{\rm sim}}[M]$ and
$\GB^{{\rm dis}}[M]$ is only in the choice of the distance function
$d$.\nl
2) Below in $\LL_{\sigma_0}$ the formula $d(x,y)=0$ appears which
normally means $x=y$. if not we could have in Def. 2.3(3), (4)
replaced $\lq k\leq r\rq$, $\lq s\leq n\rq$ by $\lq k<r\rq$, $\lq
s<n\rq$ respectively.
\lz
{\bf 2.3 Definition:} 1) For a system $\GB=(M, I, h, d):$\nl
for $x\in I$,

$N_r(x)\eqdf\{y\in I: d(y, x)\le r\}$\nl
for $x\in M\cup I,$

$ N_r^+(x)\eqdf\{ y\in M\cup I:d(h(y), h(x))\le r\}$. \nl
2) $\LL_\sigma$ is the set of first order formulas for $\sigma$; we have
variables on $M$, variables on $I$, the predicates of $\tau_1, \tau_2$; and
the additional atomic formulas $h(x)=y$; and $\lq h(y)\in N_r (h(x))\rq$ for
each $r$ (see part (4) below).  \nl
3) $\GB^-=\big( M^{\GB}, I^{\GB}, h\big)$, $\GB_r=\big( M^{\GB},
I^{\GB}, h^{\GB}, \lq d(x, y)\le k\rq\big)_{k\le r}$ so $\GB_0=\GB^-$ and
we let $\sigma_r=\sigma_r(\GB)=\sigma(\GB_r)$ (so, $\GB_r$ is a two
sorted model but not a system).
\nl
4) So $\SL_{\sigma_n}$ is defined as in part (2) but $\lq h(y)\in
N_s(h(x))\rq$ appear only for $s\le n$.
\nz
\centerline{$*$\qquad $*$\quad$*$}
We usually apply our theorems in the following case:
\lz
{\bf 2.4 Definition:} A system  $\GB$ is {\it simple} if \nl
(a) $d(x, y)\le 1$ (where $x,y$ vary on $I$) is equivalent to a quantifier
free formula in $\LL_{\sigma_0(\GB)}$.  \nl
(b) for $x, y\in I$ we have: $d(x,y)\le r $ {\it iff} there are $x_0, \dots,
x_r\in I$ such that $x=x_0$, $y=x_r$, and $d(x_\ell, x_{\ell+1})\le 1$ for
$\ell<r$ (i.e., like 2.2(5)).
\lz
{\bf 2.4A Remark:} If $\GB$ is simple, note that every formula in
$\SL_{\sigma(\GB)}$ is equivalent to some formula in $\SL_{\sigma_0}$, if we
know just that clause (b) of Definition 2.4 is satisfied then every formula
in $\SL_{\sigma(\GB)}$ is equivalent to some formula in $\SL_{\sigma_1}$.
\lz
{\bf 2.5 Convention:} 1) We define $f$: we let $f_n(r)=r+3^n$ for $r, n\in
\bbn$ \underbar{or} more generally, $f$ a two place function (written
$f_n(r)$) from $\bbN$ to $\bbN$ satisfying: $f_n$ non decreasing in $n$ and
in $r$, $r<f_n(r)\in \bbn$ and $f_n^{(3)}(r)\le f_{n+1}(r)$ where
$f_n^{(0)}(r)=r$, $f^{(\ell+1)}_n(r)=f_n(f_n^{(\ell)}(r))$ and
$f^{(2)}_n(r)\ge f_n(r)+f_n(0)$.\nz
2) We call $f$ {\it nice} if in addition $f_n^{(4)}(r) \le f_{n+1}(r)$. \nl
3) For $g$ a function from $\bbN$ to $\bbN$, let $g(\lk r_\ell:\ell< m\rk)=
\lk g(r_\ell):\ell<m\rk$.
\lz
{\bf 2.6 Definition:} 1) For a system $\GB= (M, I, h, d)$ and
$m, n\in \bbn$, and $\bar a=\langle a_0, \dots, a_{m-1}\rangle\subseteq \GB$
we define $\th^n_r(\bar a, \GB)$, here $r$ stands for a distance.
We define it by induction\footnote{$^\dagger$} {In the following
definition basic is atomic or a negation of atomic.} on $n$:
$$\eqalign{\th^0_r(\bar a, \GB)\; {{\rm is }}\; \{\varphi (x_0, \dots,
x_{m-1}): & (M, I, h)\models \varphi [a_{0},\dots, a_{m-1}], \cr & \varphi\;
{{\rm a \;basic\;formula\; in\;}} \LL_{\sigma(\GB_r)}\}\cr \th^{n+1}_r(\bar
a, \GB)= \{\th^n_r(\bar a\conc \langle c\rangle, \GB): & c\in M\cup
I\}.\cr}$$

\noindent  2) If $\sigma$ is a vocabulary of systems, $n, m\in \bbn$,
\underbar{then} $\TH^n_r(m ,\sigma)$ is the set of formally possible
$\th^n_r(\bar a, \GB)$ for $\GB$ a $\sigma$-system, $\bar{a}\subseteq\GB$
and ${\rm lg}(\bar{a})=m$; this is defined naturally. Pedantically, we define it by
induction on $n$; \underbar{for $n=0$} it is the family of sets $\bt$ of
basic formulas $\varphi(x_0, \dots, x_{m-1})$ of $\SL_{\sigma_r}$ such that
for each atomic $\varphi(x_0,\dots, x_{m-1})$ exactly one of $\varphi(x_0,
\dots, x_{m-1}),\,\neg \varphi(x_0, \dots, x_{m-1})$ belongs to $\bt$. \nl
TH$_r^{n+1}(m,\sigma)$ is the family of subsets of TH$^n_r(m+1,\sigma)$.\nl
3) If $\tau$ is a vocabulary of models (see 2.2(4)), $n, m\in \bbn$
\underbar{then} $\TH^n(m, \tau)$ is the set of formally possible
$\th^n_0(\bar a, \GB)$, $\GB$ a $\tau$-model, i.e.  $\GB=\GB\,{}^\rsim [M]$
for some $\tau$-model $M$ (note: the value of $r$ is immaterial as the
distance function is trivial).  \nl
4) If $r=0$ we may omit it, so for a model $M$, using $I=M$, $h$ the
identity we get the usual $th^n(\bar a, M)$ (but we do not assume knowledge
about it).\nl
5) If $\bar a$ is empty sequence we may omit it.
\lz
{\bf 2.7 Claim:} 1) For $\GB, n , m$, $\bar a$ as above,
$\varphi=\varphi(x_0,\dots, x_{m-1})$ a (first order) formula in
$\LL_{\sigma (\GB_r)}$ of quantifier depth $n$, we have: \nl
from $\th^n_r(\bar a , \GB)$ we can compute the truth value of
$\lq\GB_r\models \varphi[\bar a]\rq$.  \nl
Here and in later instance we mean:\nl
for any $\bt\in \TH^n_r(m, \sigma)$ we can compute a truth value
$\Gt$ such that: if $\bt=\th^n_r(\bar a, \GB)$ then $\Gt$ is the truth
value of $\lq\GB_r \models \varphi[\bar a]\rq$. Also in the proof we
behave similarly.\nl
2) For any $\sigma$ and $r,n,m\in\bbn $, if $\bt\in \TH^n_r(m, \sigma)$
\underbar{then} for some formula
$\varphi(x_0,\dots, x_{m-1})\in \LL_{\sigma_r}$ of quantifier depth $n$, for
any $\sigma$-system $\GB$ and
$\bar a=\langle a_0, \dots, a_{m-1}\rangle \subseteq \GB$ we have:
$\bt=\th^n_r(\bar a, \GB)$ iff $\GB\models \varphi[\bar a]$.  \nl
3) The functions $\lk \sigma, n, m, r\rk \mapsto \TH^n_r(m, \sigma)$ and
$\langle\tau, n, m\rangle\mapsto $TH$^n(m, \tau)$ are computable.  \nl
4) From $\th^n_r(\bar a, \GB)$ we can compute $\th^m_s(\bar a, \GB)$ if $r\ge
s$ and $n\ge m$. Also if $\Rang(\bar b)\subseteq \Rang(\bar a)$ then from $
\th^n_r(\bar a, \GB)$ and $\{(\ell, m):b_\ell=a_m\}$ we can compute
$\th^n_r(\bar b, \GB)$. (See part (1).)\nl
5) If $\GB$ is simple \underbar{then} from $\th^{n+r}_1(\bar a, \GB)$ we can
compute $\th^n_{2^r}(\bar a, \GB)$.  \nl
6) If $n_1, n_2\ge 2^d$ \underbar{then} $\th^d(\lk\rk, (n_1, <))=
\th^d(\lk\rk, (n_2, <))$.  Also if $\th^d(\lk\rk, M^\ell_i)=\bt$ for
$\ell=1, 2$ and for
$i<\max\{n_1, n_2\}$, $\tau$
%=\{<, P_0, \dots, P_{m+1}\}$ where each $P_\ell$
%is a unary predicate
a vocabulary, $M_i\in K_\tau$ {\it then} $\th^d(\lk\rk,
\sum_{i<n_1}M^1_i)= \th^d(\lk\rk, \sum_{i<n_2} M^2_i)$.
If $M^\ell_i\in K_\tau$ for $\ell=1,2$, $i<k$ and $\th^d(M^1_i)=
\th^d(M^2_i)$ for $i<k$ then $\th^d(\langle\rangle, \sum_{i<k} M^1_i) =
\th^d(\langle\rangle, \sum_{i<k}M_i^2)$.\nz
7) For a given vocabulary $\tau$, and $d\in \bbn$ there is an operation
$\oplus$ on $\TH^d(0, \tau)$ such that $\th^d(\lk\rk, \sum_{i<k}M_i)= \oplus
\lk \th^d(\lk\rk, M_i):i<k\rk$, this operation is associative (but in
general not commutative).
\hz
{\it Proof:} 1) We prove this by induction on the formula. (It goes
without saying that the reasoning below does not depend on $\GB$.)
\lz
{\ku $\varphi$ atomic:}
Thus $n=d(\varphi)=0$, and by the Definition of $\th^n_r(\bar a, \GB)$
the statement is trivial.
\lz
{\ku $\varphi=\neg\psi$:} Easy by the induction hypothesis.
\lz
{\ku $\varphi=\varphi_1\wedge\varphi_2$} (or $\varphi_1\vee\varphi_2$,
or $\varphi_1\to \varphi_2$): Easy by the induction hypothesis.
\lz
{\ku $\varphi=(\exists x)\varphi_1$:} Without loss of generality
$\varphi=(\exists x_m)\varphi_1(x_0, \dots, x_{m-1}, x_m$).
So $d(\varphi_1)=n-1$, and by the induction hypothesis for $a_0, \dots,
a_m\in \GB$ we have: the truth value of $\GB\models \varphi_1[a_0, \dots,
a_{m-1}, a_m$] is  computable from $\th^{n-1}_r(\lk a_0, \dots, a_m\rk,
\GB)$. Say it holds iff $\th^{n-1}_r(\lk a_0, \dots, a_m\rk, \GB)\in
\bT_{\varphi_1}$ ($\bT_{\varphi_1}$ a subset of
$\TH^{n-1}_r(m+1,\sigma(\GB)))$. \nz
Now $\GB\models \varphi[a_0, \dots, a_{m-1}$] iff for some $a_m$,
$\GB\models \varphi_1[a_0, \dots, a_{m-1}, a_m]$, iff \nz
$\th^{n-1}_r(\lk a_0, \dots, a_{m-1}, a_m\rk, \GB)\in \bT_{\varphi_1}$
($\bT_{\varphi_1}$ the subset of ${\rm TH}^{n-1}_r(m+1,\sigma(\GB))$
from above) for
some $a_m \in \GB$ {\ku iff} for some $c\in \GB$, $\th^{n-1}_r(\lk a_0,
\dots, a_{m-1}\rk\conc\lk c\rk, \GB)\in \bT_{\varphi_1}$ {\ku iff}
$\th^n_r(\lk a_0, \dots, a_{m-1}\rk,\GB)$ is not disjoint to
$\bT_{\varphi_1}$ ( the first ``iff" by the definition of satisfaction,
the second ``iff" by the choice of $\bT_{\varphi_1}$, the third ``iff" is
trivial; the last ``iff'' by the induction step in the definition of
$\th^n_r$). So we have completed the induction.
\nz
2) We define $\varphi=\varphi_{\bt}$ for $\bt\in \TH^n_r(m, \sigma)$ as
required by induction on $n$; check the inductive definition of
$\th^n_r$.
\nz
3) Read the definition. (2.6(2))
\nz
4) By induction on $n$, for $n=0$ as $\th^n_r$ ``speaks'' on more
basic formulas. For $n+1$ using the induction hypothesis (and the
definition of $\th^{n+1}_r$).
\nz
5) We prove it by induction on $n$. The step from $n$ to $n+1$ is very
straightforward. For $n=0$, we prove the statement by induction on $r$.
For $r=0$ note $\th^n_{2^r}(\bar a,\GB)=\th^0_{2^0}(\bar a,
\GB)=\th^0_1(\bar a, \GB)$ so there is nothing to prove. For $r=r(0)+1$ just
note that for $s_0\le 2^{r}$ we have: $d(x, y)\leq s_0$
is equivalent to $\bigvee_{{s_1+s_2=s_0\atop s_1,s_2\le 2^{r(0)}}}(\exists
z)[d(x, z)\le s_1\,\&\,d(z, y)\le s_2$]
\nz
6) The first phrase is a special case of the second
(with $M_i$ a model with a single element:
$i$). Let $n(\ell)\mathop{=}\limits^{\rm def} n_\ell$.
For the second phrase we prove the following more
general statement by induction
on $d$:

\item{$(*)_d$} Assume that for $\ell=1, 2$ we have:
$$M_\ell=\sum_{i<n(\ell)} M^\ell_i \hbox{ and }0\leq i_\ell(1)<i_\ell(2)<
\dots <i_\ell(k^*-1)<n(\ell),$$
we stipulate $i_\ell(0)=-1$ and $i_\ell(k^*)=n(\ell)$ (possibly $k^*=1$),
assume further $\bar a^\ell_k\subseteq M^\ell_{i_\ell(k)}$ has length
$m(k)$ for $k< k^*$.
Also assume that for each $k=0, \dots, k^*-1$ we have $\th^d(\bar a^1_k,
M^1_{i_1(k)})= \th^d(\bar a^2_k, M^2_{i_2(k)})$ and $(i_1(k+1)-i_1(k)-1)$ and
$(i_2(k+1)-i_2(k)-1)$ are equal or both $\ge 2^d-1$. Lastly assume
$\th^d(M^1_i)=\th^d(M^2_j)$ at least when $(\exists m)[i\in (i_1(m),
i_1(m+1))\&j\in (i_2(m), i_2(m+1)]$ (holds automatically when proving
second phrase of (6)).
\underbar{Then}
$$\th^d(\bar a^1, \sum_{i<n(1)} M^1_i)=\th(\bar a^2, \sum_{i<n(2)} M^2_i)
\ {\rm where }\ \bar a^\ell=\bar a^\ell_0\conc\bar a^\ell_1\conc\bar
a^\ell_2\conc\dots \conc \bar a^\ell_{k^*-1}.$$

\noindent
The proof is straightforward, and for the case $k^*=1$ we get the desired
conclusions. Lastly the third phrase of (6) is also a particular case of
$(*)_d$: let $n(1)=n(2)=k$, $k^*= n(\ell)+1$, $i_\ell(m)= m-1$, and $\bar
a^\ell_i=\langle\rangle$.\nz
7) The proof is like that of (6), but in $(*)_d$ we add
$\bigwedge_{k=1}^{k^*-1}i_1(k)=i_2(k)$.
\lz
{\bf Remark:} If we will want to quantify on sequence of elements (i.e.
use $\bar c$ rather than $c$) this helps.
\hfill\qed$_{2.7}$
\lz
{\bf 2.7A Claim:} Let $\GB$ be a system, let $\langle a_i:i<m\rangle$
and $\langle r_i:i<m\rangle$ (where $a_i\in \GB$ and $r_i\in \bbn)$
$f_n:\bbn\to\bbn$ for $n\in\bbn$ (not necessarily as in 2.5) be given
\underbar{then} for some $\langle n_i:i\in w\rangle$ where $w\subseteq
m$ ($=\{0, \ldots, m-1\}$)
and function $g$ from $\{0, 1, \dots, m-1\}$ to $w$ which is the identity on
$w$ we have: $\sum_{i\in w} n_i\le m-\vert w\vert$ and the sets in $\langle
N^+_{f_{n_i}(\sum\{r_j: g(j)=i\})}(a_i):i\in w\rangle$ are pairwise disjoint
and $\bar a$ is included in their union provided that
\item{$(*)$} $2f_{n_1}(r_1)+f_{n_2}(r_2)\le f_{n_1+n_2+1}(r_1+r_2)$ and
$f$ non decreasing in $r$ and in $n$ (considering $f$ a two place function
from $\bbn$ to $\bbn$).
\hz
{\ku Proof:} We call $\lk w, \bar n , g\rk$ with $\bar n=\lk n_i:i\in w\rk$
\underbar{candidate} if it satisfies all the requirements in the
conclusion of 2.7A  except possibly the ``pairwise disjoint''. Clearly there
is a candidate: $w=\{0, \dots, m-1\}$, $\bigwedge_{i\in w}n_i=0$, $g$ the identity
on $w$. So there is a candidate $\lk w, \bar n, g\rk$ with $\vert w\vert$
minimal. If the disjointness demand holds, then we are done.
So assume $i(1)\neq i(2)$ are in $w$ and there is $x$ belonging to
$N^+_{f_{n_{i(1)}}(\sum\{r_j:g(j)=i(1)\})}(a_{i(1)})$ and to
$N^+_{f_{n_{i(2)}}(\sum\{r_j:g(j)=i(2)\}}(a_{i(1)})$.

Let $w'=w\setminus \{i(2)\}$, $g'$ be a function with domain $\{0,\dots,
m-1\}$ defined by: $g'(j)$ is $g(j)$ if $g(j)\neq i(2)$, and $g'(j)$ is
$i(1)$ if $g(j)=i(2)$. Lastly define $n'_i$ for $i\in w'$: $n'_i=n_i$ if
$i\neq i(1)$ and $n'_i=n_{i(1)}+n_{i(2)}+1$ if $i=i(1)$ and let $\bar n'=\lk
n'_i:i\in w'\rk$. Now we shall show below that $(w', \bar m', g')$
is a candidate thus finishing the proof: for this we have to check the
two relevant conditions.
First
$$\sum_{i\in w'}n'_i=\sum_{{i\in w'\atop i\neq i(1)}}n'_i+n'_{i(1)}=
\sum_{{i\in w'\atop i\neq i(1)}}n_i+n_{i(1)}+n_{i(2)}+1=$$
$$\sum_{i\in w}n_i+1\le m-\vert w\vert+1=m-(\vert w\vert -1)=m-\vert
w'\vert.$$
Secondly, why $\bigcup_{i\in w} N^+_{f_{n'_i}(\sum\{r_j:g(j)=i\})}(a_i)$
includes $\bar a?$ if $j(*)<m$ then for some $i\in w$ we have $a_{j(*)}\in
N^+_{f_{n_i}(\sum\{r_j:g(j)=i\})}$; if $i\neq i(1), i(2)$ then $i\in
w'$, and
$$N^+_{f_{n'_i}(\sum\{r_j:g'(j)=i\})}(a_i)=N^+_{f_{n_i}(\sum\{r_j:g(j)=
i\})}(a_i)$$
so we are done; if $i=i(1)$ then $n'_i=n_{i(1)}+ n_{i(2)} +1 \geq
n_{i(1)} = n_i$ and
$$\sum \{r_j:g'(j)=i\}\ge\sum\{r_j:g(j)=i\}$$
hence
$$N^+_{f_{n'_i}(\sum\{r_j:g'(j)=i\})}( a_i)\supseteq
N^+_{f_{n_i}(\sum\{r_j:g(j)=i\})}(a_i)$$
and we are done.

We are left with the case $i=i(2)$, so by the choice of $i$ we have
$d(a_{j(*)}, a_{i(2)})\le f_{n_{i(2)}}(\sum\{r_j:g(j)=i(2)\})$ and by the
choice of $x$ (and $i(1), i(2))$ above $d(a_{i(2)}, x)\le
f_{n_{i(2)}}(\sum\{r_j:g(j)=i(2)\})$ and $d(x, a_{i(1)})\le
f_{n_{i(1)}}(\sum\{r_j:g(j)=i(1)\})$.  So as $d$ is a metric (i.e. the
triangular inequality) $d(a_{j(*)},a_{i(1)})\le 2
f_{n_{i(2)}}(\sum\{r_j:g(j)=i(2)\})+ f_{n_{i(1)}}(\sum\{r_j:g(j)=i(1)\})$.
Now $\sum\{r_j:g'(j)=i(1)\}=
\sum\{r_j:g(j)=i(1)\}+\sum\{r_j:g(j)=i(2)\}$ (by the definition of $g'$) and
$n'_{i(1)}=n_{i(1)}+n_{i(2)}+1$, hence what we need is
\item{$(*)$} $2 f_{n^1}(r^1)+f_{n^2}(r^2)\le f_{n^1+n^2+1}(r^1+r^2)$

\noindent which is assumed.
\hfill\qed$_{2.7A}$
\lz
{\bf 2.7B Remark:} 1) We can replace in 2.7A and its proof
$\sum\{r_j:g(j)=i\}$ by $\max\{r_j:g(j)=i\}$ and
 (in $(*)$) $r^1+r^2$ by $\max\{r^1, r^2\}$.
\nz
(2) Concerning $(*)$,
letting $n=\max\{n^1, n^2\}$ and
$r=\max \{r^1, r^2\}$ it suffices to have
\item{$(*)_1$} $f_n(r)$ is non decreasing in $n$ and in r.
\item{$(*)_2$} $3f_n(r)\le f_{n+1}(r)$
\lz
{\bf 2.8 Definition:} 1) For a system $\GB$, and $r, n, m\in \bbn$, and
$\bar a =\langle a_0, \dots, a_{m-1}\rangle\subseteq \GB$ we call $\bar a$ a
$(\GB, r)$-component \underbar{if} $\bar a\subseteq N^+_r(a_0)$.  In this
case we define $\bth^n_r(\bar a, \GB)$ ($\bth$ is for bounded theory).  We
do it by induction on $n$ (for all $r$) (the function $f$ from 2.5 is an
implicit parameter.)\nz
$(\alpha)$ $\bth^0_r(\bar a, \GB)=\th^0(\bar a,\GB_r)$ \nz
$(\beta)$ $\bth^{n+1}_r(\bar a, \GB)$ is $\langle \bt_0, \bt_1,
\bt_2\rangle$ where
$$\eqalign{
{\rm(i)}\; \bt_0&=\bth^n_r(\bar a, \GB)\cr
{\rm (ii)}\;   \bt_1&=\{\bth^n_r(c, \GB):
   c\in N^+_{f^{(2)}_n(r)}(a_0)\setminus N^+_{f_n(r)}(a_0)\}\cr
{\rm (iii)}\;   \bt_2&=\{\bth^n_{f^{(2)}_n(r)} (\bar a\conc \lk c\rk, \GB):
   c\in N^+_{f_n(r)}(a_0)\}\cr
}$$
2) If $\sigma$ is a vocabulary of systems and $r, n, m\in \bbn$ then
BTH$^n_r(m, \sigma$) is the set of formally possible $\bth^n_r(\bar a,
\GB)$, $(\GB$ a $\sigma$-system, $\bar a=\lk a_0, \dots,
a_{m-1}\rangle\subseteq \GB$ and $\bar a$ is a $(\GB, r)$-component).
\lz
{\bf 2.9 Claim:} 1) For any $\sigma$ (vocabulary of systems), numbers $n, m,
r\in \bbn$, and $\bt\in \BTH^n_r(m, \sigma)$, \underbar{there is}
$\varphi=\varphi(x_0, \dots, x_{m-1})\in \LL_\sigma$ (even
$\SL_{\sigma(\GB_{f_n(r)})}$) of quantifier depth $n$ such that for any
$\sigma$-system $\GB$, and $\bar a=\<a_0, \dots, a_{m-1}\>\subseteq \GB$ a
$(\GB, r)$-component we have:
$$\bt=\bth^n_r(\bar a,\GB) {\ \ \rm iff\ \ }
\GB_{f_n(r)}\restrictedto N^+_{f_n(r)}(a_0)\models \varphi[\bar a].$$
\nz
2) If $n\ge m$, $\bar b$ is the permutation of $\bar a$ by the
function $h$ or $\bar b=\lk a_\ell:\ell<k\rangle$ for some $k\leq
\lg(\bar{a})$  and $b_0= a_0$ 
\underbar{then}
from $\bt=\bth^n_r(\bar a, \GB)$ we can compute $\bth^m_r(\bar b, \GB)$
(using $n,m,r$, and $h$ or $k$).
\hz
{\it Proof:} Should be clear (see convention 2.5(1)).
\lz
{\bf Remark:} Concerning 2.9(2) we can say something also in the case
$b_0 \neq a_0$ but there was no real need.
\lz
{\bf 2.10 Definition:} For a $\sigma $-system $\GB$, and $r, n,m\in \bbn$ we
define $\GB^n_{r, m}$ as the expansion of $\GB$ by the relations
$R^\ell_{\bt}=\{\bar a:\bar a$ is a $(\GB, r')$-component ,
$\bt=\bth^n_{r'}(\bar a, \GB)$, $[\ell=1\Rightarrow \bar a\subseteq M]$ and
$[\ell=2\Rightarrow \bar a \subseteq I]\}$ for each
$\bt\in \BTH^n_{r'}(m', \sigma)$
$m'\leq m+n$, $r'\leq f^{(3)}_n(r)$ and
$\ell\in \{1, 2\}$.  We let $I^n_{r, m}[\GB]=\GB^n_{r,
m}\restrictedto I$.  Writing $\bar r=\langle r_\ell:\ell<k\rangle$ we mean
$\max (\bar r)$.  Writing $\GB^n_{\bar r, \bar m}$ means
the common expansion of
$\GB^n_{r_\ell, m_\ell}$ for $\ell<\lg(\bar r)=\lg(\bar m)$ if
$\lg(\bar r)=0$ we mean $\GB^n_{0,0}$ (we could
alternatively use $\GB^n_{\max(\bar r),\max(\bar m)}$, make little
difference).
Writing $\le r $ we mean for every $r'\le r$.
\lz
{\bf 2.11 Claim:} $I^{n+1}_{\bar r, \bar m} [\GB]$ essentially expands
$I^n_{\bar r', \bar m'}[\GB]$ when
\itemitem{(ii)} $\bar r' =r\conc\langle 0\rangle$, $\bar m'= \bar m\conc
\langle 1\rangle$ or
\itemitem{(ii)} $\bar r'\leq f^{(2)}_n(\bar r)$, $\bar m'\leq \langle
m_i+1: i< {\rm lg}(\bar m)$

\noindent
(essentially expand means that every predicate in the latter is equivalent
to a quantifier free formula in the former, the function giving this is
the scheme of expansion).
\hz {\it Proof:} Should be
clear by $(i)$ of $(\beta)$ of 2.8(1).
\lz
{\bf 2.12 Definition:} For a system $\GB$ and $n, m \in \bbn$ and $\bar
r=\lk r_\ell:\ell<m\rk$ such that $r_\ell\in \bbn$ and $\bar a=\langle a_0,
\dots, a_{k-1}\rangle\subseteq I^{\GB}$ \nl
1) We say $\bar a$ is $(n, \bar r)$-sparse for $\GB$ if $\bar r=\langle
r_\ell:\ell<m\rangle$ and
$$\ell<k<m\ {\rm and }\ \ N^+_{f_n(r_\ell)}(a_\ell)\cap
N^+_{f_n({r_k})}(a_k)=\emptyset$$
moreover (slightly stronger)
$$
d(a_\ell, a_k)\geq f_n(r_\ell)+ f_n(r_k)+1.
$$
2) We define $\uth^n_{\bar r}(\bar a, \GB)$ for an $(n, \bar r)$-sparse
$\bar a\subseteq\GB$, by induction on $n$:\nz
$\uth^0_{\bar r}(\bar a , \GB)=\th^0(\bar a, \GB) $ and $\uth^{n+1}_{\bar
r}(\bar a, \GB)=\langle \bt_0, \bt_1\rangle$ where:
$$
\bt_0=\{\langle \bar s, \uth^n_{\bar s}(\bar a, \GB)\rangle: \bar
s\leq f^{(2)}_n(r)\}
$$
(see 2.5, $f_n^{(2)}(\bar
r)=\lk f_n(f_n(r_\ell)):\ell<\lg (\bar a)\rk$, remember
$f^{(2)}_n(\bar r)\le
f_{n+1}(\bar r)$)
$$
\eqalign{\bt_1=\{\langle \bar s,
\uth^n_{\bar s\conc\langle 0\rangle}(\bar a\conc \lk c\rk,\GB)\rangle: &
\bar s\le \bar r\hbox{ (i.e. }\bigwedge_\ell s_\ell\le r_\ell\hbox{) and
}c\in I^{\GB}\hbox{ and }\cr
 & \bar a\conc \langle c\rangle
\hbox{ is }(n, \bar s \conc \< 0\>)\hbox{-sparse i.e. }\cr
 & N^+_{f_n(0)}(c)\hbox{ is disjoint from }N^+_{f_n(s_\ell)}
(a_\ell)\hbox{ for }\ell<\lg (\bar a) \}.\cr}
$$
3) $\UTH^n_{\bar r}(m, \sigma)$ is the set of formally possible
$\uth^n_{\bar r}(\bar a, \GB)$ $(\bar a$ is $(n, \bar r)$-sparse for $\GB$
of length $m$, $\GB$ a $\sigma$-system, etc.).
\lz
{\bf 2.13 Claim:} 1) For every $\bt\in \UTH^n_r(m, \sigma)$ there is a
formula $\varphi=\varphi(x_0,\dots,x_{k-1})\in\LL_{\sigma_{f^*_n(\bar r)}}$
where $f^*_n(\bar r)= \min\{m:$ if $\ell_1<\ell_2<\lg(\bar a)$ then
$m>f_n(r_{\ell_1})+f_n(r_{\ell_2})\}$ of quantifier depth $n$ such that:
For a $\sigma$-system $\GB$ and $(n, \bar r)$-sparse $\bar a =\langle a_0,
\dots, a_{m-1}\rangle\subseteq I^{\GB}$, we have: $\GB\models \varphi [a_0,
\dots, a_{m-1}]$ \underbar{iff} $\bt= \uth^n_{\bar r}(\bar a, \GB)$ (and
being $(n, \bar r)$-sparse is equivalent to some quantifier free formula).\nl
2) In Definition 2.12 only $\GB\restrictedto I$ matters and in part (1),
the quantifications are on $I$ only.\nl
3) If $\GB'$ essentially expands $\GB$ then from (the scheme of
expansion and) $\uth^n_{\bar r}(\bar a, \GB')$ we can compute
$\uth^n_{\bar r}(\bar a, \GB)$.
\hz
{\it Proof:} Should be clear.
\lz
{\bf 2.14 Main Lemma:} Let $\sigma$ be a system-vocabulary; if $\otimes_0$
holds then $\otimes_n$ holds for every $n$ where: \nl
$\otimes_n$ there are functions $F_{n, \bar r, \bar m}$, for $ \bar r=\lk
r_\ell:\ell<k\rk$, $\bar m=\lk m_\ell:\ell<k\rk$, where $k, r_\ell,
m_\ell\in\bbn$, such that: \nl
$(*)$ \underbar{if} $\GB$ is a $\sigma $-system, $\bar a=\bar a^0\conc\dots
\conc \bar a^{k-1}, $ $\bar a^\ell$ an $(\GB, r_\ell)$-component,
$\<a_0^\ell:\ell<k\>$ is $(n, \bar r)$-sparse and $\bar m=\<\lg (\bar
a^\ell):\ell<k\>$ \underbar{then} letting
$\bt_\ell\mathop{=}\limits^{\rm def} \bth^n_{r_\ell}(\bar
a^\ell, \GB)$ (for $\ell<k$) and $\bt \mathop{=}\limits^{\rm
def} \uth^n_{\bar r} \Big(\<h
(a_0^\ell):\ell<k\rk, I^n_{\bar r,\bar m}[\GB]\Big)$ we have $\th^n_0(\bar
a, \GB)=F_{n, \bar r, \bar m} (\bt,\bt_0, \dots, \bt_{k-1}$).  \nl
$(**)$ $F_{n, \bar r, \bar m}$ is recursive in its variables, $n,\bar r,
\bar m$ and the functions $F_{0, \bar r',\bar m'}$ where $\bar
m'=\<m'_i:i<k'\rangle$, $k'\ge k$, $\bar r'\le f^*_n(\bar r)$ 
(see 2.13(1)) and
for $i<k$ we have $m'_i\ge m_i$ and $\sum_{i<k}(m_i'-m_i)+\sum^{k'-1}_{i=k}
m_i'\le n$.
\lz
{\bf 2.14A Remark:} 
%1) 
Why we use $\th^n_0$ and not $\th^n_r$ (in the
conclusion of $(*)$)? We can if we assume
$\lq d(x, y)\le s\,\&\,x\in I\, \&y\in
I\rq$ is an atomic formula for $\GB$ for $s\le r$.\nz
%2) Instead of using $t_{\bar s}$ for all $\bar s\le \bar r$ we can use some
%or change Definition 2.12(2) to help in clause ($\delta$) in the first part
%of the proof below.
\lz
{\it Proof:} We prove this by induction on $n$, (for all $\bar r,\bar m$),
so for $n=0$ clearly $\otimes_n$ holds.  So assume $\otimes_n$ and we shall
prove $\otimes_{n+1}$.  We shall now describe a value $\bt$ computed from
$\bth^{n+1}_{r_\ell}(\bar a^\ell, \GB)$ ( for $\ell< k)$ and
$\uth^{n+1}(\langle h(a^\ell_0):\ell<k\>, I^{n+1}_{\bar r,\bar
m}[\GB])$. 
Our
intention is that $\bt=\th^{n+1}(\bar a, \GB)$.  Remember $\bt=\{\th^n(\bar
a\conc \lk c\rk, \GB):c\in \GB\}$.

Now $\bt$ will be the union of two sets.\nl
We use an informal description as it is clearer.\nl
{\bf First Part:} The set of $\th^n(\bar a\conc\< c\>, \GB$), where for some
$\ell(*)<k, c\in N^+_{f_n^{(2)}(r_{\ell(*)})} (a^{\ell(*)}_0)$.  Why can we
compute this (using the induction hypothesis of course), for each such $c$
we have: let $\bar b^\ell=\bar a^\ell$ if $\ell<k, \ell\neq \ell(*)$ and let
$\bar b^{\ell(*) }=\bar a^{(\ell(*))}\conc\<c\>$ (i.e.
$b^{\ell(*)}_{m_{\ell(*)}}=c$); next $r'_\ell$ is: $r_\ell$ for $\ell\neq
\ell(*)$ and $f^{(2)}_n (r_{\ell(*)})$ if $\ell=\ell(*)$. Necessarily for
$\ell<k$, $m^i_\ell $ is $m_\ell $ if $\ell\neq \ell(*)$ and $m_{\ell(*)}+1$
if $\ell=\ell(*)$.  Now:
\item{($\alpha$)} $N^+_{f_n(r'_\ell)} (b^\ell_0)$ for $\ell<k$ are
pairwise disjoint, moreover for $\ell(1)< \ell(2)< k$, $d(b^{\ell(1)}_0,
b^{\ell(2)}_0) > f_n(r'_{\ell(1)}) + f_n ( r'_{\ell(2)})$. \nz
[Why? for $\ell=\ell(*)$ remember that for every $r$:
$f_n(f_n^{(2)}(r))=f_n^{(3)}(r))\le f_{n+1}(r)$ and for $\ell\neq \ell(*)$
remember that for every $r$: $f_n(r)\le f_{n+1}(r)$ so in all cases clearly
$f_n(r'_\ell)\le f_{n+1}(r_\ell)$ and $b_0^\ell=a_0^\ell$.  So if
$\ell(1)\neq\ell(2)$ are $<k$, $N^+_{f_n(r^1_{\ell(1)})}(a^{\ell(1)}_0)\cap
N^+_{f_n(r'_{\ell(2)})}(b^{\ell(2)}_0)=\emptyset$ as: if
$\ell(1)\neq\ell(*)$, $\ell(2)\neq \ell(*)$, trivial.  Otherwise without
loss of generality $\ell(2)=\ell(*)$, as $\bar a$ satisfies the assumption
of $(*)$ ($(*)$ is from the lemma) $\bar a$ is $(n+1, \bar r)$-sparse
hence (see Def. 2.12(1))
$d(a^{\ell(1)}_0,a^{\ell(*)}_0)>f_{n+1}(r_{\ell(1)})+ f_{n+1}
(r_{\ell(2)})$, hence $d(b_0^{\ell(1)}, b_0^{\ell(*)})= d(a_0^{\ell(1)},
a^{\ell(*)}_0)>f_{n+1}(r_{\ell(1)})+f_{n+1}(r_{\ell(2)})\ge
f_n(r_{\ell(1)})+f_n(f_n^{(2)}(r_{\ell(*)}))=f_n(r'_{\ell(1)})+
f_n(r'_{\ell(*)}))$, as required.]
\item{($\beta$)} $\bar b^\ell\subseteq N^+_{r'_\ell}(b_0^\ell)$, so
$\bar b^\ell$ is an $r'_\ell$-component.
\nz
[Why? when $\ell=\ell(*)$ as $r_\ell\le f_n^{(2)}(r_\ell)=
r'_\ell$ and assumption on $c$; for $\ell\neq \ell(*)$ trivial.]
\item{$(\gamma)$} We can compute $\bth^n_{r'_\ell}(\bar b^\ell, \GB) $
for $\ell\neq \ell(*)$ \nz
[Why? by monotonicity properties of $\bth$ i.e. by 2.9(2) (that is
clause (i) i.e. $\bt_0$ from Def. 2.8(1)$(\beta)$)];
\item{$(\delta)$} We can compute the set of possibilities of
$\bth^n_{r'_{\ell(*)}} (\bar b^{\ell(*)}, \GB)$. \nz
[Why?  those possibilities are listed in $\bth^{n+1}$ (see $\bt_2$ in
Definition 2.8(1)($\beta)$ mainly clause(ii)).]
\item{($\varepsilon)$} We can compute $\uth^n_{\bar r'}
(\<h(b_0^\ell):\ell<k\>, I^n_{r, \bar m'}[\GB])$. \nz
 [As we can compute
$\uth^n_{\bar s} (\langle h(b^\ell_0): \ell<k\rangle, I^{n+1}_{\bar r,
\bar m}[\GB]$) for $\bar s\leq
f^{(2)}_n(\bar r)$ by the definition of $\bt_0$ in Def. 2.12(2), choose $\bar
s= \bar r'$; now, by 2.11, $I^{n+1}_{\bar r, \bar m}[\GB]$ essentially
expand 
$I^{n}_{\bar r', \bar m'}[\GB]$ 
(see clause (ii) there) hence by 2.13(3) we can get the required
object.]
\item{($\zeta$)} We can compute the set of possibilities of $\th^n(\bar
b^0\conc\ldots\conc \bar b^{k-1}, \GB)$ gotten as above for fixed
$\ell(*)<k$, all $c\in N^+_{f^{(2)}_n(r_{\ell(*)})}(a^{\ell(*)}_0)$.\nz
[Why? by $(\alpha)$ - $(\varepsilon)$ above and $\oplus_n$.]
\item{$(\xi)$} We can compute $\{\th^n(\bar a\conc c, \GB):$ for some
$\ell(*)< k$, $c\in N^+_{f^{(2)}_n(r_{\ell(*)})}(a^{\ell(*)}_0)\}$

\noindent {\bf Second Part:} The set of $\th^n (\bar a\conc\<c\>, \GB)$
where for each $\ell<k$, $c\notin N^+_{f^{(2)}_n(r)+
f_n(0)}(a^\ell_0)$.\nl
Why
can we compute this?  for such a $c$ we can let $k'=k+1$, $\bar a^k=\<c\>$
(so $c=a^k_0$, $m_k=1$).  Let $r'_\ell=r_\ell,$ (for $\ell<k$) and $r'_k=0$.
Let $\bar{m}'$ be $\bar{m}\conc\langle 1\rangle$. Now
\item{($\alpha$)} $N^+_{f_n(r_\ell)}(a^\ell_0)$ for $\ell\le k$ are
pairwise disjoint and $d(a^{\ell(1)}_0, a^{\ell(2)}_0) >
f_n(r_{\ell(1)}) + f_n( r_{\ell(2)})$.\nz
[Why? as $f_n(r)\le f_{n+1}(r)$ for $\ell(1)< \ell(2)<k$ this is
trivial, and for $\ell(1)=\ell < k=\ell(2)$ we have
$d(a^{\ell(1)}_0, a^{\ell(2)}_0)=d(a^\ell_0,
a^k_0)=d(a_0^\ell, c) > f^{(2)}_n(r_\ell)\ge f_n(r_\ell)+f_n(0)=
f_n(r'_\ell)+f_n(r'_k)$, as required].
\item{($\beta$)} $\bar a^\ell\subseteq N^+_{r'_\ell} (a^\ell_0)$, so
$\bar a^\ell$ is an $r'_\ell$-component. \nz
[Why? for $\ell<k$ as $r'_\ell=r_\ell$, for $\ell=k$ as $\bar a^k=\lk
a^k_0\rk (=\lk c\rk)]$.
\item{$(\gamma)$} we can compute $\bth^n_{r'_\ell}
(\bar a^\ell, \GB)$ for $\ell<k$ \nz
[Why? by monotonicity properties of $\bth$ i.e. 2.9(2)].
\item{($\delta$)} we can compute the possibilities for pairs $(\bt',\bt'')$
where $\bt'= \bth^{n+1}_0(\<c\>, \GB)$ and $\bt''= \uth^n_{\bar r'} (\lk
h(a_0^\ell):\ell<k\>\conc \< h(c)\>, I^{n+1}_{\bar r, \bar{m}}[\GB])$. \nz
[Why? straightforward; by the definition of $\uth^{n+1}$, i.e. $\bt_1$
of Def. 2.12(2) and $I^n_{\bar
r, \bar{m}'}[\GB]]$.\nz
\item{$(\varepsilon)$} in $(\delta)$ we can replace $I^{n+1}_{\bar r,
\bar m}[\GB]$ by $I^n_{\bar r', \bar m'}[\GB]$.\nz
[Why? by 2.11 and 2.13(3).]

\noindent
By the induction hypothesis this is enough. \nz
{\it Why the union of the two parts is $\th^n(\bar a, \GB)$?}\nl
Both obviously give subsets, and if $c$ fails the first part then
$\ell<k \Rightarrow c\notin 
N^+_{ f^{(2)}_n(r_\ell)} (a_0^\ell)$.  So
$N^+_{f_n(0)}(c)$ is disjoint to such $N^+_{{f_n(r_\ell)}}(a_0^\ell)$
and moreover $d(a^\ell_0, c)> f_n(r_\ell) + f_n(0)$.\nz
\null\hfill\qed$_{2.14}$
\lz

{\bf 2.15 Conclusion:} For any system vocabulary $\sigma$, and (first order)
sentence $\varphi$ of quantifier depth $n$, given $F_{0, \bar r,\bar m}$'s
satisfying $\otimes_0$ of 2.14 we can compute numbers $r$, $m$ and
a sentence $\psi_\varphi$ of
quantifier depth $n$, (whose vocabulary is that of $I^n_{r, m}[\GB]$)
such that: \nl
$(*)$ if $\GB$ is a $\sigma$-system which
satisfied $\otimes_0$ as exemplified by $\<F_{0, \bar r,\bar m}:\bar r, \bar
m\> $ \underbar{then} $\GB\models \varphi \Leftrightarrow I^n_{r, m}
[\GB]\models \psi_\varphi$.\hz
{\it Proof:} By 2.14 and 2.13(1) and 2.7(1) (and see 2.14A(2)).
\lz
{\bf 2.16 Conclusion:} Let $\tau$ be a vocabulary (finite for simplicity)
including a binary relation $<$, and $\varphi\in \LL_\tau$.  \underbar{Then}
we can compute an $m<\omega$, formulas $\varphi_i(x)\in \LL_\tau$ for $i<m$
with $d(\varphi_i)\le d(\varphi)$ and a sentence $\psi_\varphi\in
\LL_{\tau^*_1}$ with $d(\psi_\varphi)\le d(\varphi)$ where
$\tau^*_1\mathop{=}\limits^{\rm def}\{<\}\cup \{P_i:i<m\}$ $(m\in
\bbn$ computable from $\varphi$, 
each $P_i$ a unary predicate) satisfying the following:
\item{$(*)$} \underbar{Assume} $M$ is a finite $\tau$-model, $<^M$ a linear
order, $P\in \tau$ is unary, such that: if $R\in \tau\setminus \{<\} $ and
$\bar a=\langle a_0, \dots, a_{n(R)-1}\rk\in R^M$, then $P^M \cap [\min \bar
a, \max \bar a]_M$ has at most one member. \nl
Define $I[M]=(P^M, \dots, P_\ell\dots)_{\ell<m}$ where $P_\ell=\Big\{a\in
P^M:\varphi_\ell[a]$ is satisfied in $M\restrictedto \{x:x\le a, \vert[x,
a]_M\cap P^M\vert\le 3^{d[\varphi]}$ or $a\le x, $ $\vert[a, x]_M \cap
P^M\vert \le 3^{d[\varphi]}\}\Big\}$ \nz
($d[\varphi]$ in the quantifier depth of $\varphi$).\nz
\underbar{Then} $M\models \varphi \Leftrightarrow I[M]\models \psi_\varphi$.
\hz
{\ku Proof:} Should be clear from 2.15.
\hz
{\bf Remark:} Concerning 2.16 we can deduce it also from \S4.
\bigbreak

%
% file 463s2.tex ends  here 
%
%%% \input 463s3
%
% file 463s3.tex starts here 
%
\tb{\S3}{Proof of the Main Theorem}
\hz

{\bf 3.1 Definition:} 1) For a finite set $J\subseteq I$ let $\spr(I, J)$ be
the set of pairs $(Q^\no, Q^\yes$), where $Q^\no\subseteq I$, $Q^\yes
\subseteq I$, $Q^\no\setminus J=Q^\yes\setminus J$, $\vert
Q^\yes\smalltriangle Q^\no\vert=1$ ($A\smalltriangle B$ is the symmetric
difference). Let $\spr(I)=\spr(I, I)$.  \nl
2) For finite $J\subseteq I$ let $\mu^*(I, J)$ be the following distribution
on $\spr(I, J)$; it is enough to describe a drawing: \nl
first choose $Q^\no\subseteq I$ (all possibilities with probability
$1/2^{|I|}$) \nl
then choose $s\in J$ (all possibilities with probability $1/\vert J\vert$) \nl
finally  let
$$Q^\yes=\cases {Q^\no\cup \{s\}& {\ku if} $s\notin Q^\no$\cr Q^\no\setminus
\{s\}& {\ku if} $s\in Q^\no$\cr}$$
We write $\mu^*(I)$ for $\mu^*(I, I)$.
\lz
{\bf 3.1A Remark:} Note that the distribution $\mu^*(I, J)$ is symmetric for
$Q^\yes, Q^\no$.
\lz
{\bf 3.2 Definition:} 1) For a linear order $I$ and $J\subseteq I$ let
$\npr(I, J)=\{(Q^\no, Q^\yes): Q^\no\subseteq Q^\yes
\subseteq I $ and $Q^\no\setminus J= Q^\yes \setminus J$ and $\vert
Q^\yes\setminus Q^\no\vert=1\}$.  \nl
2) If $I$ is a set of natural numbers we use the usual order.  \nl
3) If $J=I$ we write $\npr(I)$.  \nl
4) Any $(I, J)$ is isomorphic to some $(n, J')$ so we can use such pairs
above.\nl
5) Let $\mu^{**}(I, J)$-be the distributions $\mu^*(I, J)$ restricted to
the case $Q^{{\rm no}}\subseteq Q^{{\rm yes}}$ (i.e. to $\npr(I, J)$).
\lz
{\bf 3.3 Claim:} Let $m, d\in \bbn$ be given, $\tau=\{<\}\cup\{P_i:i<m\}$,
$P_i$ a unary predicate, $K=\{M:M$ a $\tau$-model, $<^M$ a linear order$\}$.
\nl
1) For every $\vare\in \bbr^+$ for every $k$ large enough, for the
distribution $\mu^*(k)$ on $\spr(k)=\{(Q^\no, Q^\yes):$ $Q^\no,
Q^\yes\subseteq k$, $\vert Q^\yes\smalltriangle Q^\no\vert =1\}$ we have:
\underbar{if} $M_i^{\Gt}\in K$, for $i<k$, ${\Gt}\in \{\yes\; \no\}$,
and we choose $\mu^*(k)$-randomly $(Q^\no, Q^\yes)\in \spr (k)$,
\underbar{then} \nl
($*)$ the probability of $\th^d(\sum_{i<k} M_i^{\if
(i\in Q^\yes)})= \th^d (\sum_{i<k}M_i^{\if(i\in Q^\no)})$ is at least\nl
\phantom{$(*)$} $1-\vare$ ($\th^d$ is defined as in 2.6(1) considering a
model as a system by 2.2(4)).  \nl
2) Also if we first choose $ (Q^\no_u, Q^\yes_u)\in \spr \Big([[{k+1\over
2}], k)\Big)$ as above and then (possibly depending on the result)
make a decision on a choice of
$Q^\no_d=Q_d^\yes\subseteq [0, [{k+1\over 2}])$ and let $Q^\yes =
Q_d^\yes\cup Q_u^\yes$, $Q^\no=Q^\no_d\cup Q^\no_u$
\underbar{then} $(*)$ above still holds.  \nl
3) If we choose $Q^\yes\subseteq k$ such that $\vert Q^\yes\vert= [{k+1\over
2}]$ and then $Q^\no\subseteq Q^\yes$, $\vert Q^\no\vert+1=\vert Q^\yes\vert
$ (all possibilities with the same probabilities) \underbar{then} with
probability tending to $1$ with $k$ going to $\infty$ we get that $(*)$
(of 3.3(1) above) holds.
\nl
4) The parallel of (1) holds for $\mu^{**}(k)$, $\npr(k)$.

Before proving 3.3 we define and note:
\lz
{\bf 3.4 Definition:} For the $\tau, K, d$ from 3.3, let $\zeta_k$ be the
maximal real in $[0, 1]$ such that: if $M^{\Gt}_i\in K$ for $i<r$,
$\Gt\in \{\yes, \no\}$ (where $r\in \bbn$),
$J\supseteq\{i<r:\th^d(M_i^\no)\neq\th^d(M_i^\yes)\}$ has $k$ elements (and
$J\subseteq I= r= \{0, \dots, r-1\}$),
\underbar{then} ${k\over r}(1-\zeta_k)\ge
\Prob_{\mu^*(I,J)} \Big( \th^d(\sum_{i<r} M_i^{\if (i\in Q^\no)}) \neq
\th^d(\sum_{i<r} M_i^{\if(i\in Q^\yes)})\| (Q^\no, Q^\yes)\in \spr(I,J)\Big)$.
Let $\xi_k = 1-\zeta_k$.
\lz
{\bf 3.4A Observation:} 1) $\zeta_k$ is well defined; and
$\zeta_k\le\zeta_{k+1}$ \nz
2) An alternative definition of $\zeta_k$ is that it is the maximum
real in $[0, 1]$ satisfying: $(1-\zeta_k)$ is not smaller than the
relative-probability of $\th^d\Big(\sum_{i<r}M_i^{{\rm if}(i\in Q^{{\rm
no}})}\Big)\neq\th^{d}\Big(\sum_{i<r}M_i^{{\rm if}(i\in Q^{{\rm
yes}})}\Big)$ for the probability distribution $\mu^*(I)$, under the
assumption $Q^{{\rm no}}\smalltriangle Q^\yes\subseteq J$.\nz
3) Without loss of generality in 3.4, $r\le 2k+1$ (and even $r=k$).
\lz
{\it Proof:}  Note that the number of possible $\lk\th^d(M^{\Gt}_i):i,
\Gt\rk$ is finite.\nz
1) By (2). First draw $i\in J$ (equal probability) and then use $\spr(I,
J\setminus \{i\})$. Alternatively let $J=\{i_0, \ldots,
i_{\ell(*)-1}\}$, where $0\leq i_0 < \ldots < i_{\ell(*)-1}< r$, $I=\{0,
\ldots, r-1\}$ so $|J|= \ell(*)$. Let $J_{\ell(*), \ell}
\mathop{=}\limits^{\rm def} J\setminus \{i_\ell\}$ for $\ell< \ell(*)$.
Clearly
\item{$(*)_1$} 
$$\Prob_{\mu(I, J)} \Big(\th^d(\sum_{i<r} M_i^{{\rm
if}(i\in Q^{\no})})\neq \th^d(\sum_{i<r}M_i^{{\rm if}(i\in Q^{\yes})})
\| (Q^\no, Q^\yes)\in \spr(I, J)\Big) 
$$
$$
= {1 \over {\ell(*)}}\sum_{\ell<\ell(*)}
2^{-(r-1)}\vert \{Q\subseteq I\setminus \{i_\ell\}: \th^d(\sum_{i<r}
M_i^{{\rm if}(i\in Q)}) \neq \th^d(\sum_{i<r} M_i^{{\rm if}(i\in
Q\cup\{i_\ell\})})\}\vert.
$$

\noindent
Hence if $\ell(*)=k+1$ then
\item{$(*)_2$} 
$$
\Prob_{\mu(I,J)}\Big(\th^d(\sum_{i<r} M_i^{{\rm if}(i\in
Q^\no)})\neq \th^d(\sum_{i<r} (M_i^{{\rm if}(i\in Q^\yes)}) \| (Q^\no,
Q^\yes)\in \spr(I, J)\Big)
$$
$$ 
= {1\over {k+1}} \sum_{\ell\leq k}
\Prob_{\mu(I, J_{k+1,\ell})}\Big({\cal E}_\ell \| (Q^\no, Q^\yes)\in
\spr(I, J_\ell)\Big),
$$
where
$$
{\cal E}_\ell=\big\{(Q^\no, Q^\yes)\in \spr(I, J): 
\th^d(\sum_{i<r} M_i^{{\rm if}(i\in
Q^\no)}) = \th^d(\sum_{i<r} M_i^{{\rm if}(i\in Q^\yes)})\big\}. 
$$

\noindent
Now compute.
\nz
2) Should be clear.\nz
3) By addition theory i.e. 2.7(7) and by 2).
\hfill\qed$_{3.4A}$
\lz
{\bf 3.5 Definition:} Let $c\in \bbn$ be the number of members in $\TH^d(0,
\tau)$ (the set of formally possible $th^d(\langle\rangle, M)$, $M\in K,$
see 2.6(3)).  Let $k_0\in \bbn$ be such that $k_0\to (3^{d+8})^2_{c^2}$
(exists by Ramsey theorem).
\lz
{\bf 3.6 Observation:} \nl
$(*)_1$ \ \ \ $\zeta_{k_0}\ge {1\over (k_0 2^{k_0})}$.\nl
[Why? Let $r$, $M^{\Gt}_i$, and $J$, $\vert J\vert=k_0$ be given. First
draw $Q^{{\rm yes, no}}\cap (r\setminus J)$ and assume they are equal. \nz
Now (by 3.4A(2)) it is enough to prove that now the probability of the equality
i.e. of
$$\th^d\Big (\sum_{i<r}M_i^{{\rm if}(i\in Q^{{\rm
no}})}\Big)=\th^{d}\Big(\sum_{i<r}M_i^{{\rm if}(i\in Q^{{\rm yes}})})\Big)$$
is $\ge 1/( k_0 2^{k_0})$, assuming $Q^{{\rm no}} \smalltriangle Q^{{\rm
yes}}=\{j\}\subseteq J$. For $i<j$ from $J$ let $\langle \bt^\no_{i, j},
\bt_{i, j}^\yes\rangle$ be $\langle
\th^d(M_i^\no+M_{i+1}^\no+\dots+M_{j-1}^\no)$,
$\th^d(M_i^\yes+M_{i+1}^\no+M_{i+2}^\no+\dots+M_{j-1}^\no)\rangle$.  So by
the choice of $k_0$ we can find $J'\subseteq J, \vert J'\vert=3^{d+8}$ and
$ \langle \bt_0, \bt_1\rangle$ such that for $i<j$ in $J'$, we have $\langle
\bt^\no_{i, j}, \bt^\yes_{i, j}\rangle= \langle \bt_0, \bt_1\rangle$.  Let
$J'=\{i_\ell:\ell<3^{d+8}\}$.  For each $j\le 3^{d+8}$ let
$Q_i=\{i_m:m<j\}$.  By addition theory for (first order theory) linear
order, (that is 2.7(6)) for $\ell \eqdf[3^{d+8}/2]$, $\th^d(\sum_i
M_i^{\if (i\in Q_\ell)})=\th^d (\sum_i M_i^{\if (i\in Q_{\ell+1})})$.  So
the probability for equality is at least the probability of $Q^\no=Q_\ell$,
$Q^\yes=Q_{\ell+1}$ which is $\ge 1/(k_0 2^{k_0})$].
\hfill\qed$_{3.6}$
\hz
{\bf 3.7 Observation:} For every $\ell$, $k>0$ we have
$$
\xi_{k \ell}\leq \xi_{k} \big(\sum_{j\leq \ell-1}
[{{\ell-1}\choose{j}} \xi^j_k (1-\xi_k)^{\ell-1-j}]\xi_j\big).
$$
\hz
{\it Proof:} Let us be given $r$, $M^{\Gt}_i$ for $i<r$ and $J$ as in
$(*)$ so $\vert J\vert = k\ell$. Choose $\langle I_j: j< \ell\rangle$ a
partition of $I\mathop{=}\limits^{\rm def} r = \{0, \ldots, r-1\}$ to
intervals such that for $j<\ell$, $J_j\mathop{=}\limits^{\rm def} J\cap
I_j$ has exactly $k$ members. Now first draw $Q^\no_j\subseteq I_j$ for
$j< \ell$ (with equal probabilities), second draw $s_j\in J_j$ for
$j<\ell$ (with equal probabilities) and third draw 
$R^\no\subseteq \{0, \ldots, \ell-1\}$ (with equal probabilities) and
fourth draw $j(*)<\ell$ (with equal probabilities). Let $Q^\yes_j =
Q^\no_j\smalltriangle \{s_j\}$, and $R^\yes = R^\no 
\smalltriangle \{j(*)\}$,
and $Q^\no= \bigcup \limits_{j<\ell} Q^{{\rm if}(\ell\in R^\no)}_j$,
$Q^\yes = \bigcup\limits_{j<\ell} Q^{{\rm if}(\ell\in R^\yes)}_j$.

Easily $(Q^\no_j, Q^\yes_j)$ was chosen by the distribution $\mu^*(I_j,
J_j)$ and $(R^\no, R^\yes)$ was chosen by the distribution $\mu^*(\{0,
\ldots, \ell-1\})$ and $(Q^\no, Q^\yes)$ was chosen by the distribution
$\mu^*(I, J)$. Hence it is enough to prove:
\item{$(*)$} $\Prob\Big( \th^d(\sum_{i<r} M^{{\rm if}(i\in Q^\no)})
\neq \th^d(\sum_{i<r} M^{{\rm if}(i\in Q^\yes)}_i)\Big) $ \nl
$\leq \xi_k
(\sum_{j\leq \ell-1} [ {{(\ell-1)}\choose {j}}\xi^j_k
(1-\xi_k)^{\ell-1-j}] \xi_j)$.

For $j<\ell$ let $N^{\Gt}_i = \sum_{i\in I_j} M^{{\rm if}(i\in
Q^{\Gt}_j)}$ and let $p_i\in [0,1]_{\BBR}$ be $\Prob \Big(
\th^d(N^\no_j) \neq \th^d(N^\yes_j)\Big)$ and let $A=\{j:
\th^d(N^\no_j)\neq \th^d(N^\yes_j)\}$, so the events $\lq j\in A\rq$
are independent and $p_i= \Prob(i\in A) \leq \xi_k$. Now if we make
the first and second drawing only, we know $A$, $\langle N^{\Gt}_i: i<
\ell, \Gt\in\{\no, \yes\}\rangle$ and modulo this, by the definition
of $\xi_{|A|}$ we know
$$
\Prob\Big( \th^d(\sum_{i<r} M_i^{{\rm if}(i\in Q^\no)}) \neq
\th^d(\sum_{i<r} M_i^{{\rm if}(i\in Q^\yes)})\|\hbox{ after 1st and 2nd
drawing}\Big)=
$$
$$
\Prob\Big( \th^d(\sum_{j<\ell} N_j^{{\rm if}(j\in
R^\no)}) \neq \th^d(\sum_{j<\ell} N_j^{{\rm if}(j\in R^\yes)}) \|\hbox{ after
1st and 2nd drawing}\Big)\leq \xi_{|A|}.
$$

As the events $\lq j\in A\rq$ are independent we can conclude
$$
\Prob_{\mu^*(I,J)}\Big(\th^d(\sum_{i<r} M^{{\rm if}(i\in Q^\no)}_i) \neq
\th^d(\sum_{i<r} M_i^{{\rm if}(i\in Q^\yes)})\| (Q^\no, Q^\yes) \in
\spr(I, J)\Big)
$$
$$
\leq \sum_{j\leq \ell} \Prob(|A|=j) \times {j\over \ell}\times \xi_j =
\sum_{j\leq \ell} \Big(\sum_{u\subseteq \ell, |u|=j} \prod_{m\in u} p_m
\prod_{m<\ell, m\notin u} (1-p_m)\Big)\times {j\over\ell} \xi_j.
$$
Now looking at this as a function in $p_m\in [0, \xi_k]_{\BBR}$ for
$m<\ell$, for some $\langle p^*_m: m<\ell\rangle$ we get maximal
values, and as the function is linear, $p^*_m\in\{0, \xi_k\}$, and as
$0\leq \xi_j\leq \xi_{j+1}\leq 1$, necessarily $p^*_m =\xi_k$ so 
$$
\Prob_{\mu^*(I, J)}\Big( \th^d(\sum_{i<r} M_i^{{\rm if}(i\in
Q^\no)}) \neq \th^d(\sum_{i<r} M_i^{{\rm if}(i\in Q^\yes)})\| (Q^\no,
Q^\yes)\in \spr(I, J)\Big)
$$
$$ \leq \sum_{j\leq \ell} (\sum_{u\subseteq \ell, |u| = 1} \prod_{m\in
u} p^*_m \prod_{m<\ell, m\notin u} (1-p^*_m))\times{j\over \ell}\times
\xi_j =
$$
$$ = \sum_{j\leq \ell} {{n}\choose{j}} (\xi_k)^j (1-\xi_k)^{\ell-j}
({j\over \ell}) \xi_j =
$$
 $$ = \sum_{0<j\leq \ell} {{\ell}\choose{j}} (\xi_k)^j
(1-\xi_k)^{\ell-j} {j\over \ell} \xi_j = $$
$$ = \sum_{j\leq \ell-1} [{{\ell}\choose
{j+1}}{{(j+1)}\over \ell} (\xi_k)^{j+1}
(1-\xi_k)^{\ell-j-1}] \xi_j = $$
$$\xi_k \sum_{j \leq \ell-1} [{{\ell-1}\choose {j}}
(\xi_k)^j (1-\xi_k)^{(\ell-1)-j}] \xi_j.
$$
\hfill\qed$_{3.7}$
\hz
{\bf 3.8 Observation.}
\item{1)} If $\ell$, $k>0$, $j_0\leq \ell \xi_k$ then $\xi_{k
\ell}\leq \xi_k({{1+\xi_{j_0}}\over 2})$.
\item{2)} If $\xi_k\leq 1-{1\over m}$, $\ell> k/\xi_k$ then 
$\xi_{k \ell} \leq \xi_k({{1+\xi_k}\over 2})$.

\noindent
{\it Proof:} 1) As $\xi_j\leq \xi_{j+1}$ we have
$$
\xi_{k\ell}\leq \xi_k\times (\sum_{j\leq \ell-1} [ {{\ell-1}\choose
{j}} (\xi_k)^j (1-\xi_k)^{\ell-1-j}] \xi_j)=
$$
$$
= \xi_k \times \Big(\sum_{j<j_0} [ {{\ell-1}\choose {j}}
(\xi_k)^j (1-\xi_k)^{\ell-1-j}] 
\times 1 + \qquad
$$
$$
\qquad\sum_{j\in [j_0, \ell-1)} [{{(\ell-1)!}\over {j!
(\ell-j-1)!}} (\xi_k)^j (1-\xi_k) ^{\ell-1-j}] \times \xi_{j_0}\Big)=
$$
$$
\leq \xi_k ({{1+\xi_{j_0}}\over 2}).$$

2) Follows.  
\hfill\qed$_{3.8}$
\hz
{\bf 3.8A Remark.}
Using ``the binomial distribution approach normal distribution'' and
3.6, clearly we get e.g.:

for every $\varepsilon>0$, for some $\ell_\varepsilon$, for every
$\ell\geq \ell_\varepsilon$ we have 
$$
\xi_{k\ell} \leq \xi_k\xi_{(1-\varepsilon)\ell}.
$$
%%%%%%%%%%%%%%%%%%%%%%%%%%%%%%%%%%%%%%%%%%%%%%%%
\hz
{\it 3.9 Proof of 3.3:} 1) By the definition of $\zeta_k$ and Observations
3.5, 3.6 we get that $\lim\limits_{k\to\infty}\zeta_k=1$ and we can finish
easily. \nl
2) Follows by (1) (and the addition theory see 2.7, particularly 2.7(7)) \nl
3) Similar proof and not used (e.g. imitate the proof of 3.6. First choose
$j(*)$ then we have probability $\zeta_\kappa$ for equality there if the
distribution is $\mu^*(I_{j(*)}, J_{j(*)}$), but the induced distribution is
very similar to it).\nl
4) Follows very easily. For ${\rm spr}(k)$ with probability $1/2$ we are
choosing by ${\rm npr}(k)$.  \hfill\qed$_{3.2}$
\hz
{\ku 3.10 Proof of 1.4:} Let a real $\vare>0$ and a sentence $\theta\in
\SL^{\fo}_\tau$ ($\tau$ -from 1.3) be given.  We shall define $\varphi$ below
(after $(*)_1$), and let $\psi=\psi_{\varphi}$, $\tau^*_1$ (a vocabulary)
and $m\in \bbn$ and $\varphi_i(x)$ for $i<m$ be  defined  as in 2.16 (for
the $\varphi$ here). Let $d$ be the quantifier depth of $\psi_\varphi$
(i.e. $d=d[\psi_\varphi$]). Let $k^*\in \bbn$ be large enough as in 3.3(4)
(for the given $\vare$, $m$ and $\tau^*_1$). Let $k$ be $(2k^*+2)
(3^{d[\theta]}+1)$ (we could have waived the $3^{d[\theta]}+1$).
Now choose by induction on $r\le k$, $m_r\in\bbn$ such that\nl
$(*)$ (a) $0=m_0$ \nl
\phantom{(*)} (b) $m_r<m_{r+1}<\dots$ \nl
\phantom{(*)} (c) for any $n$,

$$\vare/3>\Prob_{\mu_n}\Big(M_n\models \bigvee_r \;(\exists x\le
m_r+1)(\exists y\ge m_{r+1}-1) [xRy]\| M_n\in K_n\Big).$$
[Why this is possible? We choose $m_r$ by induction on $r$.  The
probability in question is, for each fixed $r$, bounded from above by
$\sum_{i<m_r}\sum_{j>m_{r+1}}p_{j-i}$, the sum is the tail of an
(absolutely) convergent infinite sum so by increasing $m_{r+1}$ we can make it
$<\vare/2^{r+2}$, this suffices].

Next we try to draw the model $M_n$ in another way.  Let $n>m_k$ be given; let
$n^*=n+k^*+1$.  Let $J=\{m_{(3^{d[\theta]+1})i}:0<i<2k^*\}$, and
$I=\{m_i: i< 2k^* 3^{d[\theta]+1}\}$.  
We first draw
$\GA_n$, ``a drunkard model $\GA_n$ for $n$''.  Drawing $\GA_n$ means: \nz
{\it laziness case=first case} if $i<j<n^*$, $\bigvee_{r<k}[i\le m_r+1\& j\ge
m_{r+1}-1]$ then $\{i, j\}$ is non edge (no drawing).  \nl
{\it normal case=second case:} if $i<j<n^*$ are not in the first case but
$\neg(\exists m\in I)[i\le m\le j]$ then we flip a coin and get $e_{i, j}\in
\{$yes, no$\}$ with probability $p_{j-i}$ (for yes).  \nl
{\it drunkard case=third case}: $i<j<n^*$ and no previous case apply; we
make two draws. In one we get $e^1_{i,j}\in \{$ yes, no$\}$ with probability
$p_{j-i-1}$ (for yes) in the second we get $e^2_{i, j}\in \{$ yes, no$\}$
with probability $p_{j-i}$ for yes (we may stipulate $p_0=0$).
\lz
Now for every $Q\subseteq J$ we define $M_Q[\GA_n]:$ it is a model $(A^Q, <,
P, R^Q)$ where \nl
$A^Q$ is $\{0, \dots, n^*-1\}\setminus Q$ (so $\Vert M_Q[\GA_n]\Vert$ is
$n^*-|Q|$
and usually $|Q|$ will be $k^*$ or $k^*+1$)\nl
$<$ is the usual order on $I^Q$ \nl
$P=\{m_r: r<k\}$ \nl
$R^Q$ is\footnote{${}^\dagger$}{We use $i, j$ so that without loss of
generality $i'<j'$.} $ \{(i, j): \{i, j\}=\{i', j'\}, i'<j'<n^*$, and: \nl
\phantom{$R^Q$ is $ \{i, j):$} (a) $(i',j')$ fall into
the second case above and $e_{i',j'}=$ yes \nl
\phantom{$R^Q$ is $ \{i,j):$} or (b) $(i', j')$ fall into the third case say
$i\le m\le j$ and $m\in I$ \nl
\phantom{$R^Q$ is $ \{i, j):$ or (b)} ($m$ is unique by `` not first
case'') $m\in Q$, and $e^1_{i', j'}=$ yes \nl
\phantom{$R^Q$ is $ \{i, j):$} or (c) $(i', j')$ fall into the third case
say $i\le m\le j$ and $m\in I$ \nl
\phantom{$R^Q$ is $ \{i, j):$ or (b)} ($m$ is unique by `` not first
case'') $m\notin Q$ and $e^2_{i', j'}=$ yes$\}$.

We also define a model $N[\GA_n]$: \nl
the universe:   $\{0, \dots, n^*-1\}$  \nl
{\it relations:}  $<$ the usual order \nl
\phantom{\it relations:} $R= \{\langle i, j \rangle: e_{i, j}={{\rm yes}},\;
i<j\}$ \nl
\phantom{\it relations:} $R^1=\{\langle i, j\rangle: e^1_{i, j}={{\rm yes}},
\; i<j\}$\nl
\phantom{\it relations:} $R^2=\{\langle i, j\rangle :e^2_{i, j}={{\rm yes}},
\; i<j\}$ \nl
\phantom{\it relations:} $P=\{m_r:r\le k\}$ (on $\kappa$ see above, before
$(*)$)

Observe

\item{$(*)_1$} in ($N[\GA_n], Q$) we can define $M_Q[\GA_n]$ by q.f.
formulas. \nl
So for some first order $\varphi$ depending on $\theta, \tau$ but not on $n$
(promised above in the beginning of the proof):
\item{$(*)_2$} $M_Q[\GA_n]\models \theta$ \underbar{iff} $(N[\GA_n],
Q)\models \varphi$.

\noindent By 2.16 (where $I [(N[\GA_n], Q)]$ is defined in 2.16 with $M$
there standing for $(N[\GA_n], Q)$ here, so its set of elements is $P$)
and the choice of $\psi_\varphi$ we have:
\item{$(*)_3$} $(N[\GA_n], Q)\models \varphi$ \underbar{iff} $I[(N[\GA_n],
Q)]\models \psi_\varphi$.

Looking at the definition of $I[N[\GA_n], Q]$ in 2.16 without loss of
generality
\item{$(*)_4$} $I[(N[\GA_n], Q)]\models \psi_\varphi$ iff
$(I[N[\GA_n]], Q) \models \psi_\varphi$.

\noindent Let $J=J^d\cup J^u$, $J^d$ an initial segment, $J^u$ an end
segment, $\vert J^u\vert=k^*$, $\vert J^d\vert=k^*-1$.  Now we define
further drawing; let $\mu^{**}[J,J^u]$ be the distribution from 3.2 above on
$\npr (J^u)$, and choose $(Q_0^u, Q_1^u)\in \npr (J^u)$ randomly by
$\mu^{**}(J^u)$ then choose $Q_1^d\subseteq J^d$ such that $\vert
Q^d_1\vert=k^*+1-\vert Q^u_1\vert$ with equal probabilities, and let
$Q^1=Q^d_1\cup Q^u_1$, and let $Q^d_0=Q^d_1$, $Q^0=Q^d_0\cup Q^u_0$.

Now for $\ell\in \{0, 1\}$ we make a further drawing: if $i<j$ is a pair
from the first possibility (in the drawing of $\GA_n$), we flip a coin for
${}^* e^\ell_{i, j}\in \{\yes, \no\}$ with probability $p_{\vert [i,
j)\setminus Q^\ell\vert}$ for yes.

Let $M^\ell_{Q^\ell}$ be $(A^{Q^\ell}, <, R^{Q^\ell}\cup \{(i, j), (j,
i):i<j$ and: ${}^*e^\ell_{i, j}=\yes\}$) (it depends on the choice of
$\GA_n$ and on the further drawing).

Now reflecting we see
\item{$(*)_5$} for $\ell=0, 1$, the distribution
of $M^\ell_{Q^\ell}$ is the same as that of
$(K_{n+1-\ell},\mu_{n+1-\ell})$ (from Def. 1.3).

\noindent Hence

\item{$(*)_6$} for $\ell=0,1$, $\Prob(M^\ell_{Q^\ell}\models \theta)=\Prob
(M_{n+1-\ell}\models \theta \| M_{n+1-\ell}\in K_{n+1-\ell})$.

\noindent By the choice of $m_r$'s
\item{$(*)_7$} $\Prob(M^\ell_{Q^\ell}=
M_{Q^\ell}[\GA_n])$ is $\ge 1-\vare/3$.

\noindent
By 3.3 above (used above: the drawing of $(Q_0^u, Q_1^u$) was randomly by
$\mu^{**}(J^u)$).
\item{$(*)_8$} the absolute value of the differences between the following
is $\le\vare/3$:
\nz
$\Prob([N[\GA_{n}]], Q^0)\models \psi_\varphi)
\nz
\Prob([N[\GA_{n}]], Q^1)\models\psi_\varphi)$.

\noindent So for $\ell=0,1$:
\item{(a)} $\Prob (M_{n+\ell}\models \theta \| M_{n+\ell}\in K_{n+\ell})=
\Prob(M^\ell_{Q^\ell}\models \theta)$\nl
[Why? by $(*)_6$.]
\item{(b)} $\Prob (M^\ell_{Q^\ell}=M_{Q^\ell}[\GA_n])\ge 1-\vare/3$\nl
[Why? by $(*)_7$.]
\item{(c)} $M_{Q^\ell}[\GA_n]\models \theta$ iff
$(I[N[\GA_n]], Q^\ell)\models
\psi_\varphi$\nl
[Why? by $(*)_1+(*)_2+(*)_3+(*)_4$.]

\noindent By (a)+(b)+(c) it suffices to prove that the  probabilities of
$$
\lq (I[N[\GA_n]], Q^\ell)\models \psi_\varphi\rq
$$
for $\ell=0, 1$ has
difference $<\vare/3$ but this holds by $(*)_8$.
\hfill\qed$_{1.4}$

\lz

%
% file 463s3.tex ends  here 
%
%%% \input 463s4
%
% file 463s4.tex starts here 
%
\tb{\S4} {Generalized sums and Distortions:}
\nz
We try here to explain the results on \S2 as distorted generalized sums (and
the connection with generalized sums) and later the connection to models
with distance.  First we present for background the definition and
theorem of the generalized sum.
\lz
{\bf 4.1 Definition:} Let $\tau_0, \tau_1, \tau_2$ be vocabularies of
models.  For a $\tau_0$-model $I$ (serving as an index model),
$\tau_1$-models, pairwise disjoint for simplicity $M_t(t\in I)$, and function
$F$ (explained below) we say that a $\tau_2$-model $M$ is the $F$-sum of
$\lk M_t:t\in I\rk$) in symbols $M=\oplus_F\{M_t:t\in I\}$ \underbar{if:}

(a) the universe $\vert M\vert$ of $M$ is $\bigcup_{t\in I}\vert M_t\vert$
(if the $M_t$'s are not pairwise disjoint: $\{(t, a):t\in I, a\in M_t\})$
and we define $h_\alpha: M\to I$, $h(a)=t$ if $a\in M_t$ (if not disjoint
$h(\lk t, a\rk)=t$),

(b) if $t_1, \dots, t_k\in I$ are pairwise distinct, $\bar a=\bar
a_1\conc\dots\conc\bar a_k$ and $\bar a_\ell\in M_{t_\ell}$ (finite
sequences) {\it then}
$$\tp_{\qf}(\bar a, \emptyset, M)=F\Big(\tp_{\qf} (\lk t_1, \dots, t_k\rk,
\emptyset, I), \tp_\qf (\bar a_1, \emptyset, M_{t_1}),
\dots,\tp_\qf (\bar a_k, \emptyset, M_{t_k})\Big).$$
Another way to say it is:

(b)$^\prime$ if $a_0, \dots, a_{k-1}\in M$ then $\tp_\qf
(\lk a_0, \dots, a_{k-1}\rk, \emptyset, M)=$
$$F\Big(\tp_\qf(\lk h(a_\ell):\ell<k\rk),\dots,\tp_\qf(\lk a_\ell:
\ell<k, t=h(a_\ell)\rk, \emptyset, M_t),\,\dots\Big)_{t\in \{h(a_\ell):
\ell<k\}}.$$
\lz
{\bf 4.2 Remark:} 1) So the form of $F$ is implicitly defined, it suffices
to look at sequences $\bar a_0\conc\dots\conc \bar a_k$ (in clause (b)) or
$\lk a_0, \dots, a_k\rk$ of length the arity of $\tau_2$ (if it is finite)
i.e. maximum numbers of places for $P$ predicates $P\in \tau_2$. \nl
2) We can consider a generalization where the universe of $M$
and equality are defined like any other relation.
\lz
{\bf 4.3 The generalized Sum Theorem:} In the notation above if
$\bigtriangle_n$ is the set of formulas of quantifier depth $n$ 
\underbar{then} we can compute from $F$ the following function: like $F$ in
(b) (or (b)$^\prime$) replacing $\tp_\qf $ by $\tp_{\bigtriangle_n}$.
\lz
{\bf 4.4 Discussion:} So looking at a sequence $\bar a$ from $M$, to find
its [quantifier free] types we need to know two things:\nl
($\alpha$) the [quantifier free] type of its restriction to each
$h^{-1}(\{t\})= \{b\in M: h(b)=t\}$ \nl
($\beta)$ the [quantifier free] type of the sequence of parts,
$\lk h(b_\ell):\ell<k\rk$ {\it in } $I$.
\lz
{\bf 4.5 Definition:} Now $M$ is a $d$-distorted $F$-sum of $\{M_t:t\in I\}$
if

(a) $d$ is a distance function on $I$.

(b) $\vert M\vert$ is the disjoint union of $A_t$ ($t\in I$) and for each
$t$ we have a model $M_t$ with universe $\cup\{A_s: s\in I, d(s, t)\le 1\}$
and

(c)$^+$\ \ if $b_0, \dots, b_{k-1}\in M$, then
$$
\eqalign{ & \tp_\qf (\lk b_0, \dots, b_{k-1}\rk, \emptyset, M)= \cr
& F\Big(\tp_\qf(\lk h(b_0)\dots, h(b_{k-1}), \emptyset,
I),\dots, \cr
& \qquad\quad\tp_\qf(\lk b_\ell: \ell<k, d(h(b_\ell), t)\le 1),
\emptyset, M_t),\dots \Big)_{t\in \{h(b_m): m<k\}}.\cr}
$$
{\bf 4.6 Remark:} Note: by $\lk b_\ell: \ell <k, \Pr(\ell)\rk $ we mean the
function $g$ with domain $\{\ell<k: \Pr (\ell)\},$ satisfying
$g(\ell)=b_\ell)$.
\lz
{\bf 4.7 Discussion:} Our main Lemma 2.14, generalizes the generalized sum
theorem, to distorted sum but naturally the distortion ``expands'' with the
quantifier depth.
\lz
{\bf 4.8 The Distorted Sum Generalized Lemma:} In the notation above, for $m$
let $M_t^m$ be the model with \nl
{\it universe:}\ \ $\cup\{A_s:d(s, t)\le m\}\cup\{s\in I:d(s, t)\le m\}$ \nl
{\it  relations:}\ \ those of the $M_s$'s i.e. for $R\in \tau_1$, a
$k$-place predicate we let
$$Q^{t, m}_R=\{\lk s,a_1, \dots, a_k\rk:s\in I, d(s, t)\le m+1,\lk a_1, \dots,
a_k\rk \in R^{M_s}\}$$
(so $\{a_1, \dots, a_k\}\subseteq\cup\{A_s:d(s, t)\le 1\})$
$$Q_d^{t,\ell, m} =\{(s_1, s_2):d(s_1,t)\le \ell\ {\rm and}\ d(s_2, t)\le
\ell, \ d(s_1, s_2)\le m\},$$
$$Q_h=\{(s, a): a\in A_s, s\in I,d(s, t)\le m\}.$$
We define $I^{[m,n]}$ as the expansion of $I$ by:\nz
for $i\leq m$ and $\varphi\in \bigtriangle_n(\tau_0)$
$$
Q^\ell_\varphi=\{ t\in I: M^\ell_t \models \varphi\}.
$$
Now there are functions $F_n$ and a number $m(n)=3^n$ computable from $F$ (and
$n$) such that: \nz
\item{$\otimes$} for $b_0, \dots, b_{k-1}\in M$ \nl
$\tp_{\bigtriangle_n}(\lk b_0,\dots, b_{k-1}\rk,\emptyset, M)=
F_n\Big(\tp_{\bigtriangle_n}(\lk h(b_0), \dots, h(b_{k-1})\rk,\emptyset,
I^{[m(n),n]}),\dots,$\nl
\phantom{$\tp_{\bigtriangle_n}(\lk b_0,$}
$\tp_{\bigtriangle_n}(\lk b_\ell:
\ell<k,d(h(b_\ell), t)\le m(n)\rk, \emptyset,M^{m(n)}_t), \dots\Big)_{t\in
\{h(b_m): m<k\}}$.
\lz
{\bf 4.9 Discussion:} 1) Now if $d$ is trivial:
$$d(x, g)=\cases{0& $x=y$\cr
\infty& $x\not= y$\cr}$$
then $M^n_t=M_t$, and 2.8 (the disorted generalized sum Lemma) becomes
degenerated to 4.3 (the generalized sum Lemma), more exactly a variant.\nz
2) Note that 4.8 improve on the result of \S2 in $m(n)$ not depending
on $k$. We can have this improvement in \S2 and in \S5. \nz
3) To prove 4.8 given $b_0, \ldots, b_{k-1}$ and looking for $y$ with
$$
\bt
=\tp_{\bigtriangle_{n+1}}(\langle b_0, \ldots, b_{k-1}, y\rangle,
\emptyset, M),
$$
we fix $w=\{\ell<k: d(h(b_\ell), h(y))\leq m(n)\}$.\nl
E.g. if $w\neq \emptyset$ and $\ell(*) =\min (w)$,
then the relevant properties of
$y$ are expressed in the balls
$$
\{z\in M: d(h(z), h(b_\ell))\leq m(n)\}
$$
for
$\ell\in w$ and
$$
\{z\in M: d(h(z), h(y))\leq m(n)\},
$$
all included in
the ball
$$
\{z\in M: d(h(z), h(b_{\ell(*)}))\leq 3m(n)\}.
$$
The case $w=\emptyset$
is simpler.
\lz

%
% file 463s4.tex ends  here 
%
%%% \input 463s5
%
% file 463s5.tex starts here 
%
\tb{\S5} {Models with Distance}
\nz
{\bf 5.1 Discussion:} We try here to explain the results of \S2 as
concerning a model with a distance function ``weakly suitable'' for the
model and the connection to models with a distance function for the whole
vocabulary which is suitable for the model. This is another variant of the
distorted sums.
\lz
{\bf 5.2 Context:} Let $\tau$ be a fixed vocabulary.\nl
1) Let $K$ be the class of $\GA=(M, d)$, $M$ a $\tau$-model, $d$ a distance
on $M$ (i.e. a two place symmetric function from $\vert M\vert$ to
$\omega\cup \{\infty \}$, $d(x, x)=0$, satisfying the triangular
inequality) and for simplicity $d(x,y)=0 \Leftrightarrow x=y$.\nz
2) $K^{\sut}\subseteq K$ is the class of $(M, d)\in K$ such that $d$ (which
is a distance on $M$) is suitable for the model, i.e.

$\otimes_1$ $\lk a_0, \dots, a_{k-1}\rk\in R^M,
R\in \tau\ \Rightarrow\ \bigwedge_{\ell<m<k}d(a_\ell, a_m)\le 1$.\nz
3) $K^\rsim \subseteq K$ is the class of $(M, d)\in K$ which are simple, i.e

$\otimes_2$ $ d(x, y)=\Min\{n:$ there are $z_0, \dots, z_n$ such that
$x=z_0, z_n=y$ and \nz
\phantom{mm$\otimes_2 d(x, y)=\Min\{n:$} $\bigwedge_{\ell<n} \bigvee_{R\in
\tau}(\exists \bar a\in R^M) [\{z_\ell, z_{\ell+1}\}\subseteq \Rang \bar
a]\}$. \nz
4) $K^{\ws}_F$ is the class of $(M, d)\in K$
which are $F$-weakly suitable which means:

\itemitem{$\otimes_3$} for every $m_i$ {\it if} $\bar a^i=\lk a^i_0, \dots,
a^i_{n_i-1})$ for $i<k$, $a^i_\ell\in M$, $d(a^i_0, a^i_\ell)\le m_i$ and
$i<j<k\Rightarrow d(a^i_0, a^j_0)> m_i+m_j+1$ {\it then} the quantifier free
type of $\bar a^0\conc \bar a^1\conc \dots \conc \bar a^{k-1}$ is computed
by $F$ for the quantifier free types of $\bar a^0, \bar a^1, \dots, \bar
a^{k-1}$ and of $\lk a^0_0, a^1_0, \dots, a^{k-1}_0\rk$

\noindent Note: we can strengthen the demands e.g. (for $f$ as in 2.5)

\item{$(*)$} $d(a^i_0, a^i_\ell)\le r_i$, $d(a^i_0, a^j_0)>f_0(r_i)+
f_0(r_j)+1$ or at least $\neg (\exists x,y)[d(x, a_i)\le f_0(r_i)\wedge d(y,
a_j)\le f_0(r_j) \wedge d(x, y)\le 1]$.

\noindent 5) $K^\as$ is the family of $\GB=(M, d)\in K$ which are almost
simple: ``$d(x, y)\le 1$'' is defined by quantifier free formula and $d(x,
y)=\Min\{n$: we can find $z_0, \dots, z_n, x=z_0, z_n=y$,
$d(z_\ell,z_{\ell+1})\le 1\}$.
\lz
{\bf 5.3 Discussion:} For $(M, d)\in K^\ws_F$ we want to ``separate'' the
quantification to bounded ones and to distant ones.  We can either note that
it fits the context of \S 2 or repeat it.
\lz
{\bf 5.4 Definition:}
1) For $\GB=(M, d)\in K$, $x\in M$, $m<\omega$ let $N^+_m(x)=\{y\in M: d( y
,x)\le m\}$. \nz
2) We define ``$\bar a$ is a $(\GB, r)$-component'' and $\bth^n_r(\bar a,
\GB)$ as in Definition 2.8(1), and BTH$^n_r(m, \tau)$ as in 2.8(2).\nl
3) We define $\GB^2_{n, m}$ as expanding $\GB$ by the relation
$R_{\bt}=\{\bar a: \bar a$ is a $(\GB, r)$-component, $\bt=\bth^n_r(\bar a,
\GB)\}$ for $\bt\in \BTH^n_r (m, \tau)$. \nl
4) We define ``$\bar a$ is $(\GB, \bar r)$-sparse'' and $\uth^n_{\bar r}
(\bar a, \GB)$'', $\UTH^n_{\bar r}(m, \tau)$ as in 2.12.
\lz
{\bf 5.5 Theorem:} the parallel of 2.14.
\lz
%{\bf 5.6 Gaifman's Theorem:} We can now conclude Gaifman's Theorem but we
%can improve the number to $2^n$.
%\lz

%
% file 463s5.tex ends  here 
%

\centerline {\bf References}

\nl
[B] E.W. Beth, Observations metamathematiques sur les structures simplement
ordonnes, {\it Collection de Logique Mathematique}, Serie A, Fasc 5,
Paris--Louvain 1954, 29--35.
\nl
[CK] Chang and Keisler, {\bf Model theory}, North Holland Publishing
Co, 1973, 1977.
\nl
[Gf] H. Gaifman, On local and non local properties, Proc. Herbrand Symp.
{\bf Logic Colloquium, 81}, ed. J. Stern; Studies in Logic Found. of
Math vol 107 North Holland (1982) 105--135.  
\nl 
[Gu] Y. Gurevich, Monadic second order theories, Chap XIII, {\bf Model
theoretic logics}, ed. Barwise and Feferman Springer Verlag (1985); 479--506.  
\nl 
[Lu] T. Luczak, notes.
\nl 
[LcSh 435] T. Luczak and S. Shelah, Convergence in homogeneous random
graphs, {\it Random Structures \& Algorithms}, 1995, in press.
\nl 
[ShSp 304] S. Shelah J. Spencer, On zero-one laws for random graph, {\it
Journal of A.M.S} vol. 1 (1988) 97--115.
\nl
[Mo] A. Mostowski, On direct products of theories, {\it Journal of Symbolic
Logic}, vol 17 (1952) 1--31.
\nl
[Mr] L. Marcus, Minimal models of theories of one function symbol, {\it
Israel J. of Math} vol 18 (1974) 117--131.
\nl
[Sh 42] S. Shelah, The monadic theory of order,
{\it Annals of Mathematics}, vol 102 (1975), 379--419.
\nl
[Sh F-120] S. Shelah, More on 0--1 laws, unpublished notes.
\nl
[Sh 467] S. Shelah, More on random finite models, in preparation.
\nl
[Sh 548] S. Shelah, Very weak zero one law for random graphs with
order and random binary functions, {\it Random Structures \& Algorithms},
submitted.
\nl
[Sp] J. Spencer, Survey/expository paper: zero one laws with variable
probabilities, {\it Journal of Symbolic Logic} vol 58 (1993) 1--14.
\nl
[Sh 550] S. Shelah, 0--1 laws, in preparation.
\nl
[Sh 551] S. Shelah, In the random graph $G(n,p),p=n^{-a}$: if $\psi$
has probability $0(n^{-\varepsilon})$ for every $\varepsilon > 0$ then
it has probability $0(e^{-n^\varepsilon})$ for some $\varepsilon > 0$,
in preparation. 
\nl
[F] R. Fagin, Probabilities in finite models, {\it Journal of Symbolic
Logic} vol 45 (1976) 129--141
\nl
[FV] S. Feferman and R. L. Vaught, The first-order properties of
algebraic systems, {\it Fundamenta Mathematicae} vol 47 (1959) 57--103 
\nl
[GKLT] Y.V. Glebskii, D.I. Kogan, M.I. Liagonkii and V.A. Talanov,
Range and degree of reliability of formulas in restricted predicate
calculus {\it Kibernetica} vol 5 (1969) 17-27; translation of {\it
Cybernetics} vol 5 pp 142-154.

\shlhetal
\bye